\documentclass[12pt]{article}

\usepackage{amsmath}
\usepackage{amsfonts}
\usepackage{amssymb}
\input xy
\xyoption{all}
\usepackage{srcltx}

\usepackage{hyperref}

\input xy
\xyoption{all}

\newtheorem{constr}{\bf \thesection.\arabic{abz}. Construction }

\newtheorem{theo}{\bf \thesection.\arabic{abz}. Theorem }

\newtheorem{lemm}{\bf \thesection.\arabic{abz}. Lemma}

\newtheorem{defi}{\bf  \thesection.\arabic{abz}. Definition }

\newtheorem{exa}{\bf \thesection.\arabic{abz}. Example}

\newtheorem{rema}{\bf \thesection.\arabic{abz}. Remark}

\newtheorem{conj}{\bf \thesection.\arabic{abz}. Conjecture}

\newtheorem{lemmadefi}{\thesection.\arabic{abz}. Lemma-Definition}

\newtheorem{probla}{\bf  Problem 1}

\newtheorem{problb}{\bf  Problem 2}

\newcounter{abz}[section]
\newcounter{equ}[section]
\newcommand{\abz}{\refstepcounter{abz}}
\newcommand{\equ}{\refstepcounter{equ}}

\def\ad{\mathrm{ad\,}}

\def\B{\mathbb{B}}

\def\dd{\mathfrak{d}}

\def\End{\mathrm{End}}
\def\Id{\mathrm{Id}}
\def\im{\mathop{\mathrm{im}}}

\def\Fu{\mathcal{F}un}
\def\F{\mathcal{F}}

\def\G{\Gamma}
\def\g{\mathfrak{g}}
\def\gl{\mathfrak{gl}}
\def\h{\mathfrak{h}}
\def\H{\mathcal{H}}

\def\k{\mathfrak{k}}

\def\rk{\mathop{\mathrm{rank}}}

\def\so{\mathfrak{so}}

\def\W{\mathcal{W}}

\def\SB{{S^\bullet}}
\def\Nb{{N^\bullet}}

\def\K{\mathbb{K}}
\def\N{\mathbb{N}}
\def\R{\mathbb{R}}

\def\L{\mathcal{L}}
\def\P{\mathbb{P}} \def\C{\mathbb{C}}
\def\SS{\mathcal{S}}
\def\SB{{S^\bullet}}

\def\E{\mathcal{E}}

\def\al{\alpha}

\def\Ga{\Gamma}

\def\d{\partial}

\def\la{\lambda}

\def\om{\omega}

\def\qed{$\square$}

\newcommand\op[1]{\mathop{\rm #1}}

\title{Kronecker webs, Nijenhuis operators, and nonlinear PDEs}
\author{
\and
Andriy Panasyuk\\
Faculty of Mathematics and Computer Science\\
University of Warmia and Mazury\\
ul. Słoneczna 54, 10-710 Olsztyn, Poland\\
panas@matman.uwm.edu.pl
}
\date{}

\providecommand{\bysame}{\leavevmode\hbox to3em{\hrulefill}\thinspace}
\providecommand{\MR}{\relax\ifhmode\unskip\space\fi MR }
\providecommand{\href}[2]{#2}

\begin{document}

\maketitle

\begin{center}
\emph{In memory of Stanis\l aw Zakrzewski (1951--1998)}
\end{center}

\begin{abstract}
The aim of this paper is two-fold. First, a survey of the theory of Kronecker webs and their relations with bihamiltonian structures and PDEs is presented. Second, a partial solution to the problem of bisymplectic realization of a bihamiltonian structure is given. Both the goals are achieved by means of the notion of a partial Nijenhuis operator, which is studied in detail.
\end{abstract}

\tableofcontents
%
\section*{Introduction}

A seminal paper of F. Magri \cite{m} gave rise to a notion of a \emph{bihamiltonian structure}, i.e. a \emph{pair of compatible Poisson structures} $\eta_1,\eta_2$ (here compatibility means that $\eta^\la=\eta_1+\la\eta_2$ is a Poisson  structure for any $\la$), which proved to be a very effective tool in the study of integrable systems and has been developed by many authors. F. J. Turiel \cite{t5} and I. M. Gelfand and I. Zakharevich \cite{gz1}, \cite{gz3} initiated the investigation of the local structure of pairs of compatible Poisson brackets. It turns out that there are two classes of bihamiltonian structures of principally different  nature (on the level of local geometry as well as in applications to integrable systems). The bihamiltonian structures of first class  called \emph{Jordan} (cf. Section \ref{s10}) or \emph{bisymplectic} consist of pairs $\eta_1,\eta_2$ such that in the pencil $\{\eta^\la\}$ almost all members are nondegenerate Poisson structures, i.e. inverse to symplectic forms. On the contrary, in the pencils corresponding to the second class of \emph{Kronecker} bihamiltonian structures all the members are degenerate of the same rank.

It is worth mentioning that for both Jordan and Kronecker cases the classical Darboux theorem fails: in general there is no local coordinate system in which $\eta_1,\eta_2$ simultaneously have constant coefficients. In order to understand local behaviour of Kronecker bihamiltonian structures Gelfand and Zakharevich \cite{gz2} proposed a  procedure which reduces the geometry of pairs of compatible Poisson brackets to the geometry of \emph{webs}. Recall that a classical web is a finite number of foliations in general position on a smooth manifold and that the main question in the theory of webs is to describe obstructions to simultaneous local straightening of these foliations, i.e. transforming them by a local diffeomorphism to foliations of parallel plains. The reduction mentioned consists in a passage to a local base $\B$ of the lagrangian foliation $\L=\bigcap_\la\SS_\la$, where $\SS_\la$ is the symplectic foliation of $\eta^\la$ (the Kroneckerity of the pair $\eta_1,\eta_2$ guarantees that indeed the distribution $\bigcap_\la T\SS_\la$ has constant rank and, moreover, the leaves of $\L$ are lagrangian in $\SS_\la$). The induced by $\{\SS_\la\}$ one-parameter family of foliations $\{\SS'_\la\}$ on $\B$ is called a \emph{Kronecker web}. This notion  naturally generalizes the notion of a classical web and the problem of existing of the ``Darboux coordinates'' for $\eta_1,\eta_2$ can be treated in spirit of the web theory as the problem of simultaneous straightening of the foliations $\SS'_\la$. Moreover, Gelfand and Zakharevich conjectured that the Kronecker web is a complete local invariant of a Kronecker bihamiltonian structure, that is, one can reconstruct $\eta_1,\eta_2$ from $\{\SS'_\la\}$ up to a local diffeomorphism. This conjecture was  proved by Turiel \cite{t1}, \cite{t2}.

Later another side of Kronecker webs appeared in the literature: their relation with nonlinear PDEs. First discovery was made by Zakharevich \cite{z2} who found a relation of \emph{Veronese webs} in 3 dimensions (special class of Kronecker webs) with a nonlinear second order PDE called the \emph{dispersionless Hirota equation} (originally called by Zakharevich a nonlinear wave equation). The last was studied in \cite{z2} from the point of view of twistor theory. Then M. Dunajski and W. Kry\'{n}ski \cite{dunajskiKr} found relations of the Hirota equation with the so-called hyper-CR Einstein--Weyl structures previously described by the \emph{hyper-CR nonlinear PDE}. Finally, B. Kruglikov and A. Panasyuk \cite{kruglikovP} built several series of contactly nonequivalent PDEs whose solutions are in a 1--1-correspondence with Veronese webs and which include the Hirota equation as a particular case.

This paper is intended as a survey of the above mentioned and related results. However, there are two aspects which hopefully allow to treat it partially as an independent research article. First, we introduce in full generality the notion of a \emph{patial Nijenhuis operator}, which was outlined in \cite{kruglikovP}, and use it as a convenient tool for defining and working with Kronecker webs, in particular in their relations with bihamiltonian structures and nonlinear PDEs. Second, we apply the machinery of partial Nijenhuis operators to the problem of local bisymplectic realizations of Kronecker bihamiltonian structures (see Section \ref{realbi}) generalizing results by F. Petalidou \cite{petalidou}. Also Theorem \ref{th3web} is new.

Let us overview the content of the paper. In Section \ref{sss1} we discuss relations between vanishing of the Nijenhuis torsion of linear operators  and compatibility of Lie brackets. We present some examples that motivate Definition \ref{PNO} of a partial Nijenhuis operator (PNO) $N:\h\to\g$, where $\g$ is a Lie algebra and $\h \subset\g$ is a Lie subalgebra.

Section \ref{s10} is devoted to the so-called Jordan--Kronecker decomposition theorem, a classical purely linear algebraic result on the normal form of a pair of linear operators. This result is a base of classification of bihamiltonian structures (cf. Jordan and Kronecker cases discussed above) and is also permanently used in the context of the pair of operators $N,I:\h\to\g$, where $I:\h\to\g$ is the canonical inclusion.

In Section \ref{s20} we study algebraic PNOs: we observe some important consequences of the definition (Lemma \ref{20.30}), in particular, we prove that a PNO $N$ induces a  Lie algebra structure $[,]_N$ on $\h$ compatible with the initial one and that $(N+\la I)(\h)$ is a Lie subalgebra for any $\la$.   We also formulate some sufficient or necessary and sufficient conditions on a partial operator to be a PNO (Lemmas \ref{20.40}, \ref{20.60}, Remark \ref{rema}), discuss some sufficient conditions on a restriction of a Nijenhuis operator to a subalgebra to be a PNO (Lemma \ref{20.701}). All these results are of independent interest, however their main aim are geometric applications in further sections.

In Section \ref{sAl} we recall the notions of a Lie algebroid and the related notion of linear Poisson structure. This framework is very convenient for defining the geometric version of PNOs and the construction of the canonical bihamiltonian structure related to a Kronecker web (cf. Section \ref{upCon}). We also discuss compatibility of Lie algebroid structures and the corresponding linear Poisson structures.

The central notion of this article, a geometric PNO, appears in Section \ref{sGPNO}. This  is a morphism of bundles $N:T\F\to TM$, where $\F$ is some foliation on a manifold $M$, such that the corresponding mapping induced on the spaces of sections $\G(T\F)\to\G(TM)$ is an algebraic PNO (the Lie bracket being the commutator of vector fields). This notion naturally generalizes the well-known
differential geometric conditions of vanishing of the Nijenhuis torsion of a (1,1)-tensor (in our terminology this last is a PNO with $T\F=TM$). We further generalize some results of Section \ref{s20} to the geometric context. A  new aspect with respect to the purely algebraic situation is that a PNO $N$ together with the induced Lie algebra structure $[,]_N$ on $\G(T\F)$ form a Lie algebroid structure on $T\F$, which, moreover, is compatible with the canonical one (see Lemma \ref{GPNO.10}(4,5)). The fact that the image of $N+\la I$ is a subalgebra, i. e. is the tangent bundle to some foliation $\F_\la$, indicates that geometric PNOs are related to 1-parameter families of foliations such as Veronese and Kronecker webs.

These last are the main objects of Section \ref{VerKro}. We first recall the definition of a Veronese web which is  a collection $\{\F_\la\}_{\la\in\R\P^1}$ of
foliations of corank $1$ on a manifold $M$ such that the annihilating one-form $(T_xF_\la)^\bot$ sweeps a Veronese curve in $\P T^*_xM$ for any $x\in M$. We than show that there is a 1--1-correspondence between Veronese webs and PNOs $N:T\F\to TM$ of generic type with $\F_\infty=\F$ and $T\F_0={N}\,T\F$ (Theorem \ref{theo1q}). Here the genericity of type means  that there is a sole Kronekcer block in the Jordan--Kronecker decomposition of the pair of operators $N_x,I_x:T_x\F\to T_xM$ for any $x\in M$. Next we naturally generalize this result to Kronecker webs and Kronecker PNOs, the last one admitting more than one Kronecker block in the decomposition. In Remark \ref{distr} we touch the problem of ``integrability of Veronese curves of distributions'' and its generalizations.

In Section \ref{upCon} we use the relations between Lie algebroids and linear Poisson structures established in Section \ref{sAl} to construct, given a PNO $N:T\F\to TM$, the canonical bihamiltonian structure on $T^*\F$. We then specify this construction to two particular cases: a Kronecker PNO and a Jordan PNO (i.e. a PNO with $T\F=TM$, a Nijenhuis operator).

In section \ref{constrA} we discuss in full generality relations  $\{$Kronecker  webs$\}\!\xymatrix{ \ar@/^/[r] & \ar@/^/[l] } \!\{$Kronecker bihamiltonian structures$\}$, in particular the compositions of the passages $\xymatrix{ \ar@/^/[r] & }$ and $\xymatrix{ & \ar@/^/[l] }$ in different order. This includes the procedure $\xymatrix{  & \ar@/^/[l]}$ of passing to the local base of a bilagrangian foliation and reconstruction $\xymatrix{  \ar@/^/[r] & \ar@/^/[l] \bullet}$ of a bihamiltonian structure from its Kronecker web up to a local diffeomorphism.

Section \ref{realbi} is devoted to the problem of local bisymplectic realizations of a Kronecker bihamiltonian structure. More precisely it can be formulated as follows. Let $\eta_{1,2}$ be a Kronecker bihamiltonian structure on a small open set $U \subset \R^n$. Does there exist a manifold $\overline M$ with a Jordan bihamiltonian structure $\bar\eta_{1,2}$ and a surjective submersion $p:\overline M\to U$ such that $p_*\bar\eta_{1,2}=\eta_{1,2}$? If such bisymplectic realizations exist, how many nonequivalent ones there are? We show that this problem is reduced to the following problem of ``realization of a Kronecker PNO'': (1) given a Kronecker PNO $N:T\F\to TM$ does there exist a Nijenhuis operator $\overline N:TM\to TM$ such that $\overline N|_{T\F}=N$? (2) how many locally nonequivalent ones there are? The answer to question (1) is affirmative by a result of Turiel. The answer to question (2) is rather impossible in full generality in view of great range of different nonequivalent Kronecker PNOs and Nijenhuis operators. However, in the next section we give an answer to this question in the particular case of Kronecker PNOs of generic type in $3$ dimensions thus solving the problem of local bisymplectic realizations of $5$-dimensional generic bihamiltonian structures.

More precisely, in Section \ref{realz} we prove that, given a Kronecker PNO $N:T\F\to TM$ of generic type (Veronese web) on a $3$-dimensional manifold,  in a
neighborhood of every point $p\in M$ there exists an  extension of $N$ to any of normal forms of a Nijenhis operator, necessarily cyclic. Such normal forms were obtained by Turiel and Grifone--Mehdi; they are listed in Appendix. We conjecture that the same is true in any dimension: a Kronecker PNO $N:T\F\to TM$ of generic type can be extended to any normal form of a cyclic Nijenhuis operator.

In Section \ref{Hiroo} we apply a sufficient condition for the restriction $N|_{T\F}$ of a Nijenhuis operator $N:TM\to TM$ to be a PNO (Lemma \ref{llll}) to the case $M=\R^3$, the foliation $\F$ of rank $2$ and $N$ being the simplest Nijenhuis operator with constant distinct eigenvalues. As a result we get a nonlinear second order PDE on the function $f$ defining the foliation $\F$. This is the above mentioned dispersionless Hirota equation. We further prove that any Veronese web in $\R^3$ defines a solution of this equation and, vice versa, any solution defines a Veronese web. This provides a 1--1-correspondence between Veronese webs and classes of solutions with respect to a natural equivalence relation.

Section \ref{listPDE} is devoted to generalizing these results to other types of Nijenhuis operators in $\R^3$. More precisely, we get a series of pairwise contactly nonequivalent nonlinear second order PDEs on a function $f$ of three variables. For each of these equations we establish a 1--1-correspondence between classes of their solutions and Veronese webs. A crucial ingredient in this correspondence is the solution for the realization problem of a Veronese web obtained in Section \ref{realz}.

In Section \ref{high} we discuss generalizations of the results of the two preceding sections to higher dimensions. In particular, we establish a 1--1-correspondence between (classes of) solutions of a certain system of nonlinear second order PDEs and certain Kronecker webs in $4$-dimensional case.

Finally, in Section \ref{notes} we make a short overview of related bibliography.

The notion of partial Nijenhuis operator as well as the majority of related results of this paper are based on \cite{PZ} and are products of discussions with Ilya Zakharevich, to whom the author would like to express his deep gratitude. The problem of bisymplectic realization of a bihamiltonian structure was posed to the author by Stanis\l aw Zakrzewski shortly before this prominent mathematical physicist has passed away in 1998. This paper is dedicated to his memory.

%
\section{Nijenhuis operators and compatible Lie brackets}
\label{sss1}

\abz\label{aNO}
\begin{defi}
Let $(\g,[,])$ be a Lie algebra, $N\colon \g\to\g$ a linear operator. A bilinear map $T_N\colon \g\times\g\to\g$ given by
$$
T_N(x,y):=[Nx,Ny]-N[x,y]_N,\quad x,y\in\g,
$$
where
$$
[x,y]_N:=[Nx,y]+[x,Ny]-N[x,y],
$$
is called the {\em Nijenhuis torsion} of the operator $N$. One calls
$N$ {\em (algebraic) Nijenhuis} if $T_N\equiv 0$.
\end{defi}

This notion has its origin in the well known in differential geometry notion of the Nijenhuis torsion of a (1,1)-tensor $\widetilde N:TM\to TM$  on a smooth manifold $M$: if $\g$ is the Lie algebra $\G(TM)$ of the vector fields on the manifold $M$ with the usual commutator bracket and $N:\g\to\g$ is generated by the endomorphism $\widetilde N$ of the tangent bundle $TM$, the definition above in fact defines a (2,1)-tensor which coincides with the Nijenhuis torsion tensor.

\abz\label{weakN}
\begin{lemm}\cite{mks} Let $N:\g\to\g$ be a linear operator acting on a Lie algebra $(\g,[,])$.
\begin{enumerate}
\item The bracket $[,]_N$ is a Lie algebra bracket if and only if $\dd T_N=0$ (here we regard $T_N$ as a $2$-cochain on the Lie algebra $(\g,[,])$ with the coefficients in the adjoint module and $\dd$ stands for the corresponding coboundary operator).
\item Assume $\dd T_N=0$. Then the Lie bracket $[,]_N$ is automatically {\em compatible} with $[,]$, i.e., $\la_1[,]+\la_2[,]_N$ is a Lie bracket for any $\la_1,\la_2\in\K$, here $\K$ is the ground field.
\end{enumerate}
\end{lemm}

Pairs $([,]_1,[,]_2)$ of compatible (as in the lemma) Lie brackets on a vector space will be called {\em bi-Lie structures}. The families of Lie brackets $\{[,]^\la\}_{\la\in\K^2}$, $[,]^\la:=\la_1[,]_1+\la_2[,]_2$, $\la:=(\la_1,\la_2)$, generated by bi-Lie structures $([,]_1,[,]_2)$  are called {\em Lie pencils} \cite{bols'}.
In particular, any algebraic Nijenhuis operator on $(\g,[,])$ generates  a bi-Lie structure on $\g$, hence also
 a bihamiltonian structure on $\g^*$ (consisting of the corresponding Lie--Poisson structures).

The following two examples are essential in our further considerations.
\abz\label{exa1}
\begin{exa}\rm
Let $\g=\gl(n)$, $ A\in\g$ be a fixed matrix, $N:=L_A$ be the operator of left multiplication by $A$. Then $N$ is algebraic Nijenhuis, $[x,y]_N=xAy-yAx=:[x,_Ay]$ is a Lie bracket, brackets $[,],[,]_N$ are compatible.
\end{exa}

\abz\label{exa2}
\begin{exa}\rm
Let $\h=\so(n)$, $ A$ be a fixed symmetric matrix. Then $[,_A]$ is a Lie bracket on $\h$ compatible with $[,]$.
\end{exa}

In the second example we constructed the bracket $[,_A]$ ``by analogy"
with the first example. It is natural to ask whether one can include
this bracket into a framework  similar to that of Nijenhuis operators, i.e. whether $[,_A]=[,]_N$ for some $N$.  Note that for general symmetric $A$ and  $N=L_A\colon \gl(n)\to\gl(n)$ and  we have $N\so(n)\not
\subset \so(n)$. However we observe that $[x,y]_N=[x,_Ay]$ for any $x,y\in \so(n)$. In order to understand what happens, assume for a moment that  $A$ is nondegenerate, i.e., $N$ is invertible.  Although $N\so(n)\not \subset \so(n)$  the subspace $N\so(n)$ is a Lie subalgebra in $\gl(n)$. From this we conclude  that $N^{-1}[Nx,Ny]\in\so(n)$ for any $x,y\in \so(n)$, i.e. $N^{-1}[N\cdot,N \cdot]$ is a new Lie algebra bracket on $\so(n)$. On the other hand, the fact that $T_N\equiv 0$ on $\gl(n)$ implies that $N^{-1}[Nx,Ny]=[x,y]_N=[x,_Ay]$, $x,y\in \so(n)$, in particular, this new bracket is compatible with the standard one.

Let us codify these considerations in a way which allows $N$ to be not invertible.

\abz\label{PNO}
\begin{defi}\rm Let $\g$ be a Lie algebra and $\h \subset \g$  a Lie subalgebra. We say that a pair $(\h,N)$, where $N\colon \h\to\g$ is a linear  operator,  is an (algebraic) {\em partial Nijenhuis operator on $\g$} (PNO for short) if the following two conditions hold:
\begin{enumerate}
\item[(i)] $[x,y]_N\in\h$ for any $x,y\in\h$;
\item[(ii)] $T_N(x,y)=0$ for any $x,y\in\h$.
\end{enumerate}
(Here $[,]_N$ and $T_N$ are  given by the same formulas as above; note that it follows from condition (i) that the  term $N[x,y]_N$ which appears in the definition of $T_N$ is correctly defined.)
\end{defi}

The examples above give the following two instances of PNOs: (1) let $\g=\h=\gl(n)$, $N=L_A$, then $(\h,N)$ is a PNO on $\g$ (which in fact is a Nijenhuis operator since $\h=\g$); (2) let $\h=\so(n),\g=\gl(n)$, $N=L_A|_\h$, where $A$ is a symmetric matrix. Then $(\h,N)$ is a PNO on $\g$.

In these examples, given a PNO $(\h,N)$ on a Lie
algebra $\g$, we obtained a  bi-Lie structure $([,],[,]_N)$ on
$\h$.  It turns out that it is also true in general (see Lemma \ref{20.30}).
Vice versa, given a bi-Lie structure
$([,],[,]_1)$ on a vector space $\h$ such that $(\h,[,])$ is a semisimple Lie algebra, one can identify $\h$ with $\tilde\h=\ad(\h) \subset\End(\h)$ and define $N:\tilde\h\to\End(\h)$ by $N:=\ad_1\circ\ad^{-1}$, where $\ad,\ad_1:\h\to\End(\h)$ are the corresponding adjoint representations. Then $N$ is a PNO; this fact was helpful in an approach to the problem of classification of
bi-Lie structures $([,],[,]_1)$ on semisimple Lie algebras
$(\h,[,])$ \cite{pBi-Lie}.

In order to study PNOs and their relations to bihamiltonian structures and Kronecker webs we  recall a classical result on normal forms of a pair of linear operators.

%
\section{The Jordan--Kronecker decomposition of a pair of linear operators}
\label{s10}

\abz\label{J--K}
\begin{theo}\cite{gantmacher}
 Consider a pair of operators $S_1,S_2\colon  V\to W$ between finite-dimensional vector spaces  over $\C$. Then there are direct decompositions $V=\bigoplus_{m=1}^nV_m$, $
 W=\bigoplus_{m=1}^nW_m$, $ S_1=\bigoplus_{m=1}^nS_{1,m}$, $S_2=\bigoplus_{m=1}^nS_{2,m}$, where $S_{j,m}\colon V_m\to W_m$, $ j=1,2$, $ m=1,\ldots,n$, such that  each 4-tuple $(S_{1,m},S_{2,m},V_m,W_m)$ is from the following list:
\begin{enumerate}
\item {\rm [the Jordan block $\mathbf{j}_\la(j_m)$]:} $\dim V_m=\dim W_m=j_m$ and in an appropriate
bases of $V_m$ and $W_m$ the matrix of $S_{1,m}$ is equal to
 $I_{j_m}$ (the unity $j_m\times j_m$-matrix) and the matrix of $S_{2,m}$ is equal to
$J_{j_m}^\la$ (the Jordan $j_m\times j_m$-block with the eigenvalue $\la$);
\item {\rm [the Jordan block $\mathbf{j}_\infty(j_m)$]:} $\dim V_m=\dim W_m=j_m$ and in an appropriate
bases of $V_m$ and $W_m$ the matrix of $S_{1,m}$ is equal to
 $J_{j_m}^0$  and the matrix of $S_{2,m}$ is equal to
$I_{j_m}$;
\item {\rm [the Kronecker block $\mathbf{k}_+(k_m)$]:} $\dim V_m=k_m$, $\dim W_m=k_m+1$ and in an appropriate
bases of $V_m,W_m$ the matrices of $S_{1,m},S_{2,m}$ are  equal to
$$
\left(\begin{array}{ccccc}
1&0&\ldots&0\\
0&1&\ldots&0\\
0&0&\ldots&0\\
 & &\ldots&   \\
0&0&\ldots& 1 \\
0&0&\ldots& 0
\end{array} \right),\qquad\left(\begin{array}{ccccc}
0&0&\ldots&0\\
1&0&\ldots&0\\
0&1&\ldots&0\\
 & &\ldots& \\
0&0&\ldots& 0\\
0&0&\ldots&1
\end{array}\right),
$$
respectively ($(k_m+1)\times k_m$-matrices);
\item {\rm [the Kronecker block $\mathbf{k}_-(k_m)$]:} $\dim V_m=k_m+1$, $\dim W_m=k_m$ and in an appropriate
bases of $V_m,W_m$ the matrices of $S_{1,m},S_{2,m}$ are  equal to
$$
\left(\begin{array}{cccccc}
1&0&0&\ldots&0&0\\
0&1&0&\ldots&0&0\\
 & & &\ldots& & \\
0&0&0&\ldots&1&0
\end{array} \right),\qquad\left(\begin{array}{cccccc}
0&1&0&\ldots&0&0\\
0&0&1&\ldots&0&0\\
 & & &\ldots& & \\
0&0&0&\ldots&0&1
\end{array}\right),$$
respectively ($k_m\times (k_m+1)$-matrices).
\end{enumerate}
\end{theo}

\abz\label{10.10}
\begin{defi}\rm The decomposition from the theorem above will be called the {\em Jordan--Kronecker} (J--K for short) decomposition of the pair $S_1,S_2$. We will call the Kronecker blocks $\mathbf{k}_+(k_m)$ ($\mathbf{k}_-(k_m)$) {\em increasing} (respectively {\em decreasing}).
\end{defi}

\abz\label{10.20}
\begin{defi}\rm Consider the pencil of operators $ \SB=\{S^\la\}$, $S^\la:=\la_1S_1+\la_2S_2$,
  $\la:=(\la_1,\la_2)$, generated by the operators $S_1,S_2\colon V\to
  W$. The set $E_\SB:=\{\la\in\C^2\mid \rk S^\la<\max_{\mu}\rk
  S^\mu\}$ will be called {\em exceptional} for $\SB$.
\end{defi}
 It is clear from the theorem above that the exceptional set $E_\SB$  is either $\{0\}$ (Kronecker case: the Jordan blocks are absent) or a finite union of lines in $\C^2$.

%
\section{Partial Nijenhuis operators (algebraic version)}
\label{s20}

In this section we consider vector spaces defined over a field $\K$ equal to $\R$ or $\C$. We study elementary properties of PNOs.

\abz\label{20.10}
\begin{defi}\rm Let $W$ be a vector space, $V \subset W$ its subspace, and $S\colon V\to W$ a linear operator. We say that a pair $(V,S)$ is a {\em partial operator} on $W$. The subspace $V$ is called the {\em domain} of $S$.
\end{defi}

Recall that algebraic PNOs were introduced in Definition \ref{PNO}.

\abz\label{20.30}
\begin{lemm}
If $(\h,N)$ is a PNO on $\g$, then:
\begin{enumerate}
\item $N\h$ is a Lie subalgebra in $\g$;
\item $(\h,N^\la)$, $N^\la:=\la_1I+\la_2N$, is a partial Nijenhuis operator on $\g$ for any $\la:=(\la_1,\la_2)\in \K^2$,  here $I\colon \h\to\g$ is the natural embedding;
\item $N^\la\h$ is a Lie subalgebra in $\g$ for any $\la$;
\item $[,]_{N^\la}$ is a Lie algebra structure on $\h$ and $N^\la\colon \h\to\g$ is a homomorphism between Lie algebras $(\h,[,]_{N^\la})$ and $(\g,[,])$;
\item the Lie bracket $[,]_N$ is compatible with the Lie bracket $[,]$ (see Lemma \ref{weakN} for the definition).
\end{enumerate}
\end{lemm}
Indeed, Item 1 is obvious. Item 2 is due to the equality $[,]_{\la_1I+\la_2N}=\la_1[,]+\la_2[,]_N$ and to the equality $T_{\la_1I+\la_2N}=\la_2^2T_N$. Item 3 follows from Items 1 and 2.

Now Items 4 and 5 follow easily from the equality $[x,y]_{\la_1I+\la_2 N}=(\la_1I+\la_2N)^{-1}[(\la_1I+\la_2N)x,(I+\la_2N)y]$, which makes sense for $(\la_1,\la_2)\not \in E_\SB$ (see Definition \ref{10.20}), where $\SB$ is the pencil of operators generated by $I,N$. \qed

\smallskip

In the following lemma we give some sufficient conditions for a partial  operator $(\h,N)$ on $\g$ to be a PNO.

\abz\label{20.40}
\begin{lemm} Let $\g$ be a Lie algebra and $\h \subset \g$ be  a Lie subalgebra. Let  $N\colon \h\to\g$ be an operator such that $N\h$ is also a Lie subalgebra. Then, if there exist $(a_k,b_k)$, $k=1,\ldots,K$, not proportional to $(1,0)$ and to $(0,1)$ such that $\h_k:=(a_kI+b_kN)\h$ is a Lie subalgebra and $\bigcap_{k=1}^K\h_k=\{0\}$, the pair $(\h,N)$ is a PNO.
\end{lemm}
For such $(a_k,b_k)$, $a_k\not=0$; put $\rho_k=b_k/a_k$. By the assumption, for any $x,y\in\h$ there exists $s=s(x,y)\in\h$ such that $[Nx,Ny]=Ns(x,y)$.  Thus
\begin{align}
[x+\rho_k Nx,y+\rho_k Ny] &= [x,y]+\rho_k([Nx,y]+[x,Ny])+\rho_k^2[Nx,Ny] \notag\\
&=(I+\rho_kN)[x,y]+\rho_k[x,y]_N+\rho_k^2Ns(x,y) \notag\\
&=(I+\rho_kN)([x,y]+\rho_ks(x,y))+\rho_k([x,y]_N-s(x,y)).\notag
\end{align}
Therefore $[x,y]_N-s(x,y)\in\h_k$ for any $k$ (since $\h_k$ is a subalgebra); hence $[x,y]_N-s(x,y)\in\bigcap_{k=1}^K\h_k=\{0\}$ and $[x,y]_N=s(x,y)\in\h$.

Now $T_N(x,y)$ of Definition \ref{aNO} is correctly defined and $T_N(x,y)=[Nx,Ny]-N[x,y]_N=Ns(x,y)-Ns(x,y)=0$.  \qed

\abz\label{20.50}
\begin{rema}\rm The idea of this lemma and its proof  is borrowed from   \cite[Theorem 4.1]{bouetoudufour}.
\end{rema}

\abz
\begin{rema}\rm Note that the assumption of existence of $(a_k,b_k)$, $k=1,\ldots,K$, such that $\h_k$ are subalgebras and
$\bigcap_{k=1}^K\h_k=\{0\}$ is a sufficient but not necessary condition for the Nijenhuis property of $N$. Say, if $N$ is a ``usual"
(i.e., $\h=\g$) nondiagonalizable\footnote{For instance, the operator of left multiplication by a nilpotent matrix on $\g=\gl(n)$.} Nijenhuis operator, then this condition is not satisfied.
Below we study for which cases the condition mentioned is also necessary (see Remark \ref{rema}) and give another necessary and sufficient conditions for the Nijenhuis property of $N$ in terms of the ``affinization" $\g[\al]$ of $\g$.
\end{rema}

\abz
\begin{lemmadefi} Let $(V,N)$ be a partial operator on a finite-dimensional
  vector space $W$ over $\C$ and let $I\colon V\to W$ be the natural
  embedding. Consider the pencil $\{ N^\la\}$, $ N^\la:=\la_1I+\la_2N$, $\la:=(\la_1,\la_2)$, generated by the operators $I,N$.  Then
\begin{enumerate}
\item The subspace $V_J:=\bigcap_{\la\in\C^2\setminus E_\Nb}\im N^\la$  lies in $V$ and is invariant w.r.t.~ $N$ (the operator $N_J:=N|_{V_J}$  will be called the {\em Jordan part} of $(V,N)$).
\item the intersection $\bigcap_{\la\in\C^2\setminus \{(0,0)\}}\im N^\la \subset V_J$, is equal to the zero subspace if and only if the Jordan part $N_J$ is diagonalizable.
\end{enumerate}
\end{lemmadefi}
Since $I$ is injective, there are no decreasing Kronecker blocks in the corresponding J--K decomposition (see Theorem \ref{J--K}). The rest of the proof is an easy consequence of the structure of this decomposition. \qed

\abz\label{rema}
\begin{rema}\rm
Now we see that the sufficient condition of ``existence of $(a_k,b_k)$, $ k=1,\ldots,K$, such that $\h_k$ are subalgebras and $\bigcap_{k=1}^K\h_k=\{0\}$"  from Lemma \ref{20.40} is necessary for the Nijenhuis property of the partial operator $(\h,N)$ on $\g$ if and only if the Jordan part $N_J$ is diagonalizable.
\end{rema}

\abz\label{20.60}
\begin{lemm} Let $\g$ be a Lie algebra and $\h \subset \g$ be  a Lie subalgebra. We write $\g[\al ]$ for the Lie algebra of polynomials with coefficients from $\g$ with the natural Lie bracket. Then a partial operator $(\h,N)$ on $\g$ is a PNO if and only if the image of the operator $N':=(I+\al N)|_{\h+\al \h}\colon (\h+\al \h)\to\g[\al ]$ is a Lie subalgebra.
\end{lemm}

Indeed, $\im N'$ is a Lie subalgebra if and only if  for any $x,y\in\h$ there exists $u=u_0+\al u_1\in \h+\al \h$ such that $[x+\al Nx,y+\al Ny]=u+\al Nu$. The left hand side of this equality is equal to $[x,y]+\al ([Nx,y]+[y,Nx])+\al ^2[Nx,Ny]$. Comparing the coefficients of different powers of $\al $ in the equality above we conclude that  $\im N'$ is a Lie subalgebra if and only if $u_0=[x,y]$, $ Nu_1=[Nx,Ny]$ and  $u_1+Nu_0=[Nx,y]+[y,Nx]$. The last three equalities are equivalent to conditions (i), (ii) of Definition \ref{PNO}. \qed

We conclude this section by studying relations between partial Nijenhuis operators and Nijenhuis operators.

\abz\label{20.701}
\begin{lemm} Let $\g$ be a Lie algebra, $\h \subset \g$  a Lie subalgebra, and $N:\g\to\g$ a Nijenhuis operator (see Definition \ref{aNO}). Assume that for some $\la\in \K$ the following two conditions hold:
(1) $\k:=(N+\la \Id_\g)\h$ is a Lie subalgebra; (2) $(N+\la \Id_\g)^{-1}(\k)=\h$ (for instance, this condition holds if $-\la$ is not an eigenvalue of $N$).

 Then $(\h,N|_\h)$ is a partial Nijenhuis operator on $\g$.
\end{lemm}
Put $N':=N+\la\Id_\g$. Due to the condition $T_{N'}=T_N\equiv 0$, for any $x,y \in\h$ we have $N'[x,y]_{N'}=[N'x,N'y]$, the last expression being an element of $\k$ by assumption (1). Hence, $[x,y]_{N'}=[x,y]_N+\la[x,y]\in\h$ by assumption (2) and also $[x,y]_N\in\h$. On the other hand, obviously $T_N\equiv 0\Longrightarrow T_{N|_\h}\equiv 0$. \qed

A natural question occurs: is it true that any partial Nijenhuis operator $(\h,N)$ on $\g$ with  $\h \varsubsetneq \g$ can be extended to a Nijenhuis operator on $\g$? We will come back to this question in Section \ref{sGPNO}.

%
\section{Lie algebroids and linear Poisson structures}
\label{sAl}

In this section we recall some notions related to Lie algebroids and linear Poisson structures, which will be used for defining the geometric version of PNOs and establishing their connections with bihamiltonian structures.

\abz\label{Al.10}
\begin{defi}\rm A {\em Lie algebroid} is
a vector bundle $E\to M$ endowed with a bundle morphism (called {\em anchor}) $\rho\colon E\to TM$ and a Lie algebra structure $[,]_E$ on the space of sections $\G(E)$ satisfying
\begin{enumerate}
\item[(i)] The induced mapping $\rho\colon \G(E)\to\G(TM)$ is a Lie algebra homomorphism (the space of vector fields $\G(TM)$ is endowed with the standard bracket; we use the same letter for the morphism of bundles and the morphism of spaces of sections).
\item[(ii)] $[x,fy]_E=f[x,y]_E+(\rho(x)f)y$ for any $x,y\in\G(E)$, $f\in\Fu (M)$ (here $\Fu(M)$ denotes the space of functions on $M$ in the corresponding category).
\end{enumerate}
\end{defi}

\abz\label{Al.20}
\begin{exa}\rm If $M=\{*\}$, then $\rho$ is trivial, $\Fu(M)=\K$ (the corresponding ground field), $\G(E)=E=\g$ is a Lie algebra.
\end{exa}

\abz\label{Al.30}
\begin{exa}\rm Let $E=TM$, $[,]_E$ be the commutator of vector fields, $\rho=\Id$. We say that $E$ is the {\em  tangent Lie algebroid} on $M$.
\end{exa}

\abz\label{Al.40}
\begin{exa}\rm Let $\F$ be a foliation on $M$. Put $E=T\F$ (the space of elements of $TM$ tangent to $\F$), $\rho=I\colon E\to TM$ for the natural inclusion, $[,]_E$ for the commutator of vector fields tangent to $\F$. We will call this Lie algebroid  structure {\em canonical}.
\end{exa}

 Given a Lie algebroid $(E,\rho,[,]_E)$, one can build a Poisson structure on $E^*$ which will be linear in fibers, i.e., the Poisson bracket $\{,\}$ of two sections of $E$ interpreted as (fiberwise) linear functions on $E^*$ will be a linear function on $E^*$ (see \cite{wcs}). If $x_1,\ldots,x_n$ are local coordinates on $M$ and $e_1,\ldots,e_r$ local basis of sections of $E$ and the corresponding structure functions are defined by
$$
\rho(e_i)=b_{ij}\frac{\partial}{\partial x_j},\qquad [e_k,e_l]_E=c_{kl}^me_m,
$$
then the linear Poisson bracket on $E^*$ is defined as
\begin{equation}\equ\label{linp}
\{x_i,x_j\}=0,\qquad \{\xi_k,\xi_l\}=c_{kl}^m\xi_m,\qquad \{\xi_i,x_j\}=-b_{ij}.
\end{equation}
Globaly, we have the following properties \cite{marle}:
 \begin{enumerate}
 \item[(1)] $\{\overline X,\overline Y\}=\overline{[X,Y]_E}$ (here $\overline X$ stands for the linear function on $E^*$ corresponding to $X\in\Ga(E)$);
     \item[(2)] $\{\overline X,q^*f\}=q^*(\rho(X)f)$ (here $q$ denotes the projection $E^*\to M$); \item[(3)] $\{q^*f,q^*g\}=0$.
\end{enumerate}
Folrmulas (\ref{linp}) show that these properties completely characterize the Poisson bracket; in other words, the Poisson bracket is completely characterized by its values on {\em linear} and  {\em base} functions.

  One can show that in fact the notions of a Lie algebroid on $E$ and of a linear Poisson structure on $E^*$ are equaivalent, i.e. they uniquely determine each other.

In the context of Examples \ref{Al.20}--\ref{Al.40}  the corresponding linear Poisson structure on $E^*$ is, respectively:
\begin{enumerate}
\item    the Lie--Poisson structure on $\g^*$;
\item the canonical nondegenerate Poisson structure $\eta_{T^*M}$ on $T^*M$;
\item   the canonical Poisson structure $\eta_{T^*\F}$ on $T^*\F$ (which is
  degenerate if dimension of leaves of $\F$ is strictly less than dimension of $M$);
   recall that $T^*\F$ is fibered into
  symplectic manifolds $T^*L$, where $L$ runs over leaves of $\F$.
\end{enumerate}

\abz\label{remanc}
\begin{rema}\rm
Note that the Poisson structure $\eta_{T^*\F}$ is completely determined by the anchor $I:T\F\to TM$ and the canonical Poisson structure $\eta_{T^*M}$; more precisely, $\eta_{T^*\F}=I^t_*\eta_{T^*M}$, where $I^t:T^*M\to T^*\F$ is the transposed map to $I$ understood as a smooth surjective submersion. Indeed, first notice that for any $X\in\Ga(T\F)$ we have the following equality of linear functions on $T^*M$: $\overline{IX}=(I^t)^*\overline X$, where $(I^t)^*$ stands for the pullback. Denote the Poisson brackets corresponding to $\eta_{T^*\F}$ and $\eta_{T^*M}$ by $\{,\}'$ and $\{,\}$ correspondingly and write $\sigma:T^*\F\to M$ and $\pi:T^*M\to M$ for the canonical projections. Then for any $X,Y\in \Ga(T\F)$ and any functions $f,g$ on $M$ we have
\begin{align*}
(I^t)^*\{\overline X, \overline Y\}'=(I^t)^*\overline{[X,Y]}=\overline{I[X,Y]}=\{\overline{IX},\overline{IY}\}=\{(I^t)^*\overline{X},
(I^t)^*\overline{Y}\}\\
(I^t)^*\{\overline X,\sigma^* f\}'=(I^t)^*\sigma^*(IXf)=\pi^*(IXf)=\{\overline{IX},\pi^*f\}=\{(I^t)^*\overline{X},(I^t)^*\sigma^*f\}\\
(I^t)^*\{\sigma^*f,\sigma^* g\}'=0=\{\pi^*f,\pi^* g\}=\{(I^t)^*\sigma^*f,(I^t)^*\sigma^* g\},
\end{align*}
which proves the claim (cf. properties (1)--(3) above).
\end{rema}

\abz\label{lIII.1}
\begin{defi}\rm Let $E\to M$ be a vector bundle with two Lie algebroid structures $([,]_1,\rho_1)$ and $([,]_2,\rho_2)$. They are called {\em compatible} if $(\la_1[,]_1+\la_2[,]_2,\la_1\rho_1+\la_2\rho_2)$ is a Lie algebroid structure for any constants $\la_1,\la_2$. Given two compatible Lie algebroid structures $([,]_1,\rho_1)$ and $([,]_2,\rho_2)$ on $E$, the family $\{(\la_1[,]_1+\la_2[,]_2,\la_1\rho_1+\la_2\rho_2)\}$ is a {\em pencil} of Lie algebroid structures on $E$.
\end{defi}

\abz\label{compat}
\begin{lemm} Let $E\to M$ be a vector bundle with two compatible Lie algebroid structures $([,]_1,\rho_1)$ and $([,]_2,\rho_2)$. Then the corresponding linear Poisson structures on the total space of $E^*$ are also compatible.
\end{lemm}
The proof easily follows from the definition of compatible algebroids and properties (1)--(3) which completely characterize the linear Poisson structure. \qed

\smallskip

One can also proceed in the other direction:

\abz\label{D2.2}
\begin{exa}\rm Let  $E=T^*M$. Assume $S_i\colon E\to TM$, $i=1,2$ are two compatible Poisson structures on $M$. Put $[x,y]_i:=\L_{S_ix}y-\L_{S_iy}x+d\langle S_ix,y\rangle,x,y\in\G(T^*M)$, $i=1,2$, for the corresponding Lie algebra structures on $\G(T^*M)$ \cite{mks}. Then $([,]_1,S_1),([,]_2,S_2)$ are compatible Lie algebroid structures on $T^*M$.
\end{exa}

However, note that these two constructions are {\em not} inverse to
each other.  Starting with Lie algebroid structures on $E$, one gets
Poisson structures on the total space $ \E$ of $E^*$.  The second construction would give Lie
algebroid structures on $T^*\E$.

%
\section{Partial Nijenhuis operators (geometric version)}
\label{sGPNO}

\abz\label{PNOdefi}
\begin{defi}\rm
Let $E=T\F$ for some foliation $\F$ on $M$.
We say that a pair $(E,N)$, where $N\colon E\to TM$ is a bundle morphism, is a (geometric) {\em   partial Nijenhuis operator} (PNO for short) on $M$ if the following two conditions hold:
\begin{enumerate}
\item[(i)] $[x,y]_N:=[Nx,y]+[x,Ny]-N[x,y]\in \G(E)$ for any $x,y\in \G(E)$ (here $[,]$ stands for the commutator of vector fields on $M$);
\item[(ii)] $T_N(x,y):=[Nx,Ny]-N[x,y]_N=0$ for any $x,y\in\G(E)$ (it follows from condition (i) that the second term is correctly defined).
\end{enumerate}
\end{defi}
In other words, a bundle morphism $N\colon E\to TM$ is a  geometric  PNO if it is an algebraic PNO regarded as a map of Lie algebras $\G(E)\to\G(TM)$ (which will be denoted by the same letter).

\abz\label{cgm}
\begin{rema}\rm This notion is very natural and probably existed in the literature earlier with no special name.
A similar notion appeared in \cite{cgm2} under the name "outer Nijenhuis tensor".

 F. J. Turiel used equivalent notion in \cite{t8,t10,t11}  in different terms. Namely, he considered a foliation $\F$ on a manifold $M$ and a morphism $N:T\F\to TM$ such that
 \begin{enumerate}\item[(1)]
 $N^*\al$ is closed along the leaves of $\F$ for any closed 1-form $\al$ satisfying $\ker\al\supset T\F$ .
 \end{enumerate}
 Then he proved that, given any extension $\overline{N}$ of $N$ to a morphism from $TM$ to $TM$,  the restriction of $T_{\overline{N}}$ to $T\F$ does not depend on the extension. So one can require that
 \begin{enumerate}\item[(2)]  $T_{\overline{N}}|_{T\F\times T\F}=0$. \end{enumerate}
We claim that in fact the two notions are equivalent, i.e. the following equivalences hold: $(i)\Longleftrightarrow(1)$, and, under the assumption that $(i)$ or $(1)$ is satisfied, $(ii)\Longleftrightarrow(2)$.
Indeed, assume that condition $(1)$ is satisfied. If $\al$ is a 1-form such that $d\al=0,\al|_{T\F}=0$, then for any vector fields $x,y$ we have $\al([x,y])=x\al(y)-y\al(x)$ and for $X,Y\in\G(T\F)$ we have $(N^*\al)([X,Y])=X(N^*\al)(Y)-Y(N^*\al)(X)$, i.e. $\al(N[X,Y])=X\al(NY)-Y\al(NX)$. Thus for any such $1$-form we have
\begin{align*}
\al([X,Y]_N)=\al([NX,Y]+[X,NY]-N[X,Y])=NX\al(Y)-Y\al(NX)+X\al(NY)-NY\al(X)\\-
X\al(NY)+Y\al(NX)=0.
\end{align*}
This implies $[X,Y]_N\in\G(T\F)$, hence condition $(i)$. These considerations are reversible and $(i)\Longleftrightarrow(1)$. Now if one of these equivalent conditions hold, $T_{\overline{N}}(X,Y)$ coincides with the expression $T_N(X,Y)$ from condition $(ii)$, is independent of the prolongation $\overline{N}$, and, obviously, $(ii)\Longleftrightarrow(2)$.
\end{rema}

Recall that the bundle $E=T\F$ has the canonical Lie algebroid structure with the canonical inclusion  $I:E \to TM$ as the anchor and the commutator of vector fields tangent to $\F$ as the Lie bracket on $\G(E)$.

\abz\label{GPNO.10}
\begin{lemm}
Let   $(E,N)$, $ N\colon E\to TM$, be a PNO on $M$. Then:
\begin{enumerate}
\item $N\G(E)$ is a Lie subalgebra in $\G(TM)$;
\item $N^\la:=\la_1I+\la_2 N$ is partial Nijenhuis for any $\la:=(\la_1,\la_2)$;
\item $N^\la\G(E)$ is a Lie subalgebra in $\G(TM)$ for any $\la$; in particular if rank of the distribution $N^\la E$ is constant, it is tangent to some foliation $\F^\la$;
\item $[,]_N$ is a Lie algebra structure on $\G(E)$ which together with the anchor $N\colon E\to TM$ form a Lie algebroid structure on $E$;
\item this new Lie algebroid structure on $E$ is compatible with the canonical Lie algebroid structure on $E$, i.e.  the family $\{([,]_{N^\la},N^\la)\}$ is a pencil of Lie algebroid structures on $E$ (see Definition \ref{lIII.1}).
\end{enumerate}
\end{lemm}
Items 1, 2, 3 are proven as in the algebraic case (Lemma \ref{20.30}). Let
us prove that $(E,[,]_{N^\la},N^\la)$ is a Lie algebroid for any
$\la$. The fact that $[,]_{N^\la}$ is a Lie algebra and that $N^\la$
is a homomorphism of Lie algebras is also proven as in algebraic
case. It remains to check the condition of compatibility of the
bracket with the anchor; by linearity it is enough to prove it with $N$ instead of $\la_1I+\la_2 N$:
\begin{align}
[x,fy]_N &=[Nx,fy]+[x,Nfy]-N[x,fy]\notag\\
  &=f[Nx,y]+((Nx)f)y+[x,fNy]-N(f[x,y]+(xf)y)\notag\\
  &=f[Nx,y]+((Nx)f)y+f[x,Ny]+(xf)Ny-N(f[x,y]+(xf)y)\notag\\
  &=f[x,y]_N+((Nx)f)y;\notag
\end{align}
note that we used only the linearity of $N$. \qed

\smallskip

The proofs of the following two lemmas follow from the corresponding lemmas in the algebraic case
(see Lemmas \ref{20.40} and  \ref{20.701}).

\abz\label{l}
\begin{lemm} Let $\F$ be a foliation on $M$. Let  $N\colon T\F\to TM$ be a vector bundle morphism such that $NT\F$ is the tangent bundle to some foliation. Then, if there exist $(\la_1^{(k)},\la_2^{(k)})$, $ k=1,\ldots,K$, linearly independent with $(1,0)$ and with $(0,1)$ such that $(\la_1^{(k)}I+\la_2^{(k)}N)T\F=T\F^{(k)}$ for some foliation $\F^{(k)}$ and $\bigcap_{k=1}^KT_x\F^{(k)}=\{0\}$ for any $x\in M$, the pair $(T\F,N)$ is a PNO.
\end{lemm}

\abz\label{llll}
\begin{lemm} Let $\F$ be a foliation on $M$. Let  $N\colon TM\to TM$ be a Nijenhuis  (1,1)-tensor such that for some $\la\in \K$ the following two conditions hold: (1) the distribution $B:=(N+\la\Id_{TM})T\F$ is tangent to some foliation; (2) $(N+\la\Id_{TM})^{-1}(B)=T\F$. Then the pair $(T\F,N|_{T\F})$ is a PNO.
\end{lemm}

\abz\label{Example 1}
\begin{exa}\rm
Let $N$ is a ``usual" Nijenhuis operator ((1,1)-tensor).  Then $N\colon E\to TM$ be a PNO with $E=TM$.
\end{exa}

Now we provide a simplest nontrivial example of a  partial Nijenhuis operator.

\abz\label{Example 3}
\begin{exa}\rm Let $M$ be any manifold and let $v,w\in \G(TM)$ be linearly independent (at each point) vector fields.
Put $E:=\langle v\rangle$, $ N\colon  E\to TM$, $ Nv:=w$. Since $E$ is a vector bundle with one-dimensional fibers, the integrability condition
on $NE$ is trivial, and it is easy to check that the operator $N$ is partial Nijenhuis.

Assume that $v,w$ are generic.  It is clear that there is no
coordinate system in which $N$ is translation-invariant.  For example,
if $ \dim M > 2$, then $ E $ and $ NE $ are not simultaneously tangent
to any 2-dimensional foliation.
\end{exa}

The 1-parameter family of foliations of rank 1 appearing in this example via Lemma \ref{GPNO.10} is an example of the so-called Kronecker web. In more details this notion is considered in the next section.

%
\section{Veronese and Kronecker webs and PNOs}
\label{VerKro}

Recall the definition of a Veronese web  \cite{gz2}.

\abz\label{Veronese}
\begin{defi}\rm Let $\{\F_s\}_{s\in\R\P^1}$ be a collection of
foliations of rank $n$ on a manifold $M^{n+1}$ of dimension $n+1$ such that in a neighbourhood of any point there exists a local coframe $\al_0,\ldots,\al_n$ with
$T\F_s=\langle \al_0+s\al_1+\cdots+s^n\al_n\rangle^\bot$ (here $\langle\cdot\rangle^\bot$ stands for the annihilator of the span
$\langle\cdot\rangle$) for any $s\in\R\P^1=\R\cup\{\infty\}$ (by definition $T\F_\infty:=\langle \al_n\rangle^\bot$).
Thus the map $\R\P^1\ni t\mapsto \langle (\al_0+s\al_1+\cdots+s^n\al_n)|_x\rangle\in\P T_x^*M$ parametrizes a {\em Veronese curve} for any $x\in M$. The whole collection $\{\F_s\}_{s\in\R\P^1}$ is a \emph{Veronese web}.
\end{defi}

It turns outs that there exists a 1-1-correspondence between Veronese webs and special PNOs. Let us say that a PNO $(T\F,{N})$ on a manifold $M^{n+1}$ is of \emph{of generic type} if the pair of operators ${N},I:T\F\to TM$, where $I:T\F\hookrightarrow TM$ is the canonical inclusion, has a unique Kronecker block
in the J--K decomposition (see Section \ref{s10}), i.e. there exist local frames $v_1,\ldots,v_n\in\Ga(T\F)$, $w_0,\ldots,w_n\in\Ga(TM)$, in which
\begin{equation}\equ\label{kronform}
I=\left[
    \begin{array}{cccc}
      1 & & &  \\
      0 & 1 &  &\\
      & \ddots&\ddots& \\
       & & 0 & 1 \\
       & & & 0 \\
    \end{array}
  \right],
{N}=\left[
    \begin{array}{cccc}
      0 & & &  \\
      1 & 0 &  &\\
            & \ddots&\ddots& \\
       & & 1 & 0 \\
       & & & 1 \\
    \end{array}
  \right].
\end{equation}

\abz\label{theo1q}
\begin{theo} There exists a 1-1-correspondence between Veronese webs $\{\F_s\}$ on $M^{n+1}$ and PNOs $(T\F,{N})$ of generic type
such that $\F_\infty=\F$ and $T\F_0={N}\,T\F$.
\end{theo}
Let $\{\F_s\}_{s\in\R\P^1}$ be a Veronese web on $M^{n+1}$. It turns out that $\{\F_s\}$ is determined by the
foliation $\F_\infty$ and an (everywhere defined) Nijenhuis operator which is built as
follows \cite{bouetoudufour,t1,t5}. Fix $s_0,\ldots,s_n\in\R$ to be
pairwise distinct nonzero numbers; for $i=0,\dots,n$ define a rank-$1$ foliation $\SS_i$ by $T_x\SS_i:=\bigcap_{j=0,j\not=i}^nT_x\F_{t_j}$, $x\in M$. Then  $T_x\SS_i+T_x\SS_k$ is an integrable distribution for any $i,k$, hence putting $\overline{N}|_{T_x\SS_i}:=s_i\Id_{T_x\SS_i}$ we will get a Nijenhuis operator.

It is easy to see that $T_x\F_{s_i}=(\overline{N}-s_iI)T_x\F_\infty$, $i=0,\ldots,n$, where $I:=\Id_{TM}$ (indeed $\ker(\overline{N}-s_iI)=T\SS_i$ is transversal to $T\F_{\infty}$ and $\im(\overline{N}-s_iI)=\sum_{j\not=i}T\SS_j=T\F_{s_i}$).
On the other hand, one can see that the map $\R\P^1\ni s\mapsto((\overline{N}-sI)T_x\F_\infty)^\bot\in \P T_x^*M$ is a Veronese curve (a priori different from the initial one).     These two curves pass through $n+2$ distinct points of $\P T_x^*M$: $n+1$ mentioned above and  $\infty$ (since $T_x\F_{\infty}=\lim_{s\to\infty}(\overline{N}-sI)T_x\F_\infty$). We conclude by the uniqueness property of the Veronese curve (Lagrange interpolation theorem) that they coincide. Hence $T_x\F_s=(\overline{N}-sI)T_x\F_\infty$ for any $s\in\R\P^1$ and $x\in M$.

By Lemma \ref{llll} (put $\la=0$) this gives us a partial Nijenhuis operator ${N}=\overline{N}|_{T\F_\infty}\colon  T\F_\infty\to TM$. Alternatively one can use Lemma \ref{l} since $\bigcap_{i=0}^nT_x\F_{s_i}=\{0\}$.

The constructed PNO $(\F_\infty,{N})$ is independent of the choice of the numbers $s_i$. Indeed, let
$(T\F_s)^\bot=\langle \al_0+s \al_1+\dots+s^n\al_n\rangle=:\langle\al^s\rangle$ and let $X_0,\dots,X_n$ be the frame
dual to the coframe $\al_0,\dots,\al_n$. Then the partial operator ${N}: T\F_\infty=\langle X_0,\dots,X_{n-1}\rangle\to TM$
satisfying $\al^s(({N}-sI)\,T\F_\infty)=0$ for any $s$ (now $I:T\F_\infty\hookrightarrow TM$ is the canonical inclusion) is uniquely determined by ${N}X_k=X_{k+1}$, $0\le k<n$.
Note also that the pair $({N},I)$ has canonical  matrix form (\ref{kronform}) in the frames $X_0,\dots,X_{n-1}$ and $X_0,\dots,X_n$.

Vice versa, let $(T\F,N)$ be a PNO of generic type on $M$. Then it is easy to see that $( N-sI)T\F=\langle \al_0+s \al_1+\dots+s^n\al_n\rangle^\bot$, where $\al_0,\dots,\al_n$ is the coframe dual to $w_0,\ldots,w_n\in\Ga(TM)$ (see (\ref{kronform})). The integrability of the distribution $(N-sI)T\F$ follows from Lemma \ref{GPNO.10}(3).
\qed

\abz\label{remareal}
\begin{rema}\rm The proof above shows that for any PNO of generic type $(T\F,{N})$ on a manifold $M$ there exists a Nijenhuis operator $\overline{N}:TM\to TM$ such that $N=\overline{N}|_{T\F}$. It turns out that such a Nijenhuis operator is not unique. The related problem of realization of PNOs of generic type is considered in Section \ref{realz} (which in turn is related to the problem of bisymplectic realizations of bihamiltonian structures, see Section \ref{realbi}).

\end{rema}

Veronese webs are particular cases of a more general notion of a
  Kronecker web \cite{z1}. Notice that F. J. Turiel (and initially the author \cite{p1}) uses the term Veronese web for both the notions \cite{t2},\cite{t8}.
\abz\label{kronweb}
\begin{defi}\rm \cite{z1}
   Let
  $\{\F_s\}_{s\in\R\P^1}$ be a collection of foliations on a manifold
  $M$. Assume that there is a vector bundle $\Phi\to M$ and two bundle
  morphisms $\phi_i\colon T^*M\to \Phi$, $i=1,2$,  such that for any
  $s_1,s_2\in\R$, $(s_1,s_2)\ne 0$ we have $\ker\phi_{(s_1,s_2)}=(T\F_{s_1\colon
  s_2})^\bot$, here $\phi_{(s_1,s_2)}:=s_1\phi_1+s_2\phi_2$. We say
  that $\{\F_s\}_{s\in\R\P^1}$ is a {\em Kronecker web} if for any
  $(s_1,s_2)\in\C^2\setminus\{(0,0)\}$ the morphism
  $s_1\phi_1+s_2\phi_2\colon (T^*M)\otimes\C\to \Phi\otimes\C$ is fiberwise
  surjective, or in other words, $\dim \ker(s_1\phi_1+s_2\phi_2)$ does
  not depend on $(s_1,s_2)\in\C^2\setminus\{(0,0)\}$ for any fixed
  point of $M$. Equivalently, the J--K decomposition of the pair of operators $\phi_{1,x},\phi_{2,x}\colon T^*_xM\to\Phi_x$, $ x\in M$, does not contain Jordan blocks (this explains the  name ``Kronecker web").
\end{defi}

It turns out  that the dualization of this definition gives an example of a
PNO.  Indeed, given a Kronecker web $\{\F_s\}_{s\in\R\P^1}$,
consider the pencil of the transposed morphisms
$\phi_{(s_1,s_2)}^t\colon \Phi^*\to TM$ (which are fiberwise injective for any
$s_{1,2}$). Note that,
$\im\phi_{(s_1,s_2)}^t=(\ker\phi_{(s_1,s_2)})^\bot=T\F_{s_1\colon
  s_2}$, in particular $\im\phi_1^t= T\F_\infty$.  Hence $\phi_1^t$
identifies $\Phi^*$ with $T\F_\infty$. Consider the map $(\phi_1^t)^{-1}\colon T\F_\infty\to \Phi^*$ and the map $N:=\phi_2^t\circ(\phi_1^t)^{-1}\colon T\F_\infty\to TM$.

We claim that $(T\F_\infty,N)$ is a PNO. Indeed, $(s_1I+s_2N)T\F_\infty=\im \phi_{(s_1,s_2)}^t= T\F_{s_1\colon s_2}$ for any $s_1\colon s_2 \in \R\P^1$, where $I$ is the canonical embedding $T\F_\infty\hookrightarrow TM$.
Moreover, one can find a finite number (which depends on the structure of Kronecker blocks in the J--K decomposition) of points in $\R\P^1$ such that the intersection of the corresponding foliations is trivial.
By Lemma \ref{l} we conclude that  $(T\F_\infty,N)$ is a PNO.

\abz\label{kronPNO}
\begin{rema}\rm One can immediately see that a Kronecker web is the
  same as a PNO $N$ such that the morphism $N^\la$ is injective at any point of the base manifold and for any $\la \ne 0$ (provided one can take complex $\la$). We will call such a PNO {\em Kronecker}, since for such $N$ the pair of morphisms $(N,I)$ contains only (increasing) Kronecker blocks in the Jordan--Kronecker decomposition at any point. Veronese webs are distinguished by the case of a sole Kronecker block (PNOs of generic type).
\end{rema}

\abz\label{realizA}
\begin{rema}\rm The proof of Theorem \ref{theo1q} suggests a question: is it true that any Kronecker PNO is a restriction to the tangent bundle of some foliation of some ``usual'' Nijenhuis operator  on $M$ as it is for the particular case of Veronese webs, see Remark \ref{remareal}.
The answer to this question is positive \cite[Theorem 2.1]{t8} (see also Sections \ref{realbi}--\ref{realz} for the discussion of the realization problem).
\end{rema}

\abz\label{distr}
\begin{rema}\rm In the context of Veronese webs the following theorem is true \cite{p2,bouetoudufour}. Let
$\al_0,\ldots,\al_n$ be a local coframe on $\R^{n+1}$ and let $D_s:=\langle \al^s\rangle^\bot$, where
$\al^s:=\al_0+s\al_1+\cdots+s^n\al_n$. Assume that the distribution of hyperplanes $D_s \subset T\R^{n+1}$ is integrable for $n+3$ different values of $s\in\R\P^1$. Then $D_s$ is integrable for any $s$, i.e., induces a Veronese web.

Note that this statement is surprising starting from  $n=3$ since the condition of integrability $d\al^s\wedge\al^s=0$ is   polynomial in $s$ of degree $2n$, thus one would expect that a sufficient condition would be vanishing of the polynomial at $2n+1$ different points.

In \cite{p2,bouetoudufour} also a generalization of this theorem was proven, considering Kronecker webs with Kronecker blocks of equal dimension.

The construction of PNO related to Kronecker webs and Lemma \ref{l} allow to prove an analogue\footnote{With $n+3$ values, where $n+1$ is the dimension of the target space of the highest Kronecker block.}
of this theorem for the most general Kronecker webs without any restrictions on the dimensions of the Kronecker blocks (another proof of such a theorem is obtained by F. J. Turiel \cite[Corollary 2.1.2]{t8}).
\end{rema}

%
\section{Canonical bihamiltonian structure related with a PNO}
\label{upCon}

Combining the construction of a linear Poisson structure from a Lie algebroid described in Section \ref{sAl} with Lemmas
\ref{GPNO.10}(4-5) and \ref{compat} one obtains, given a PNO $(T\F,N)$ on a manifold $M$, a canonically defined pencil of
(linear) Poisson structures on the total space of $T^*\F$.  We will say that this bihamiltonian structure is obtained by means of \emph{``up construction''} from a  PNO $(T\F,N)$.

Let us consider this bihamiltonian structure in detail. One of the linear Poisson structures from this pencil, $\eta_{T^*\F}$, corresponds to the canonical Lie algebroid structure on $T\F$ with the anchor $I:T\F\to TM$ (the canonical inclusion).  We know (see Remark \ref{remanc}) that $\eta_{T^*\F}=(I^t)_*\eta_{T^*M}$, where $\eta_{T^*M}$ is the canonical Poisson structure on $T^*M$. Analogous statement is true for the second generator of this pencil.

\abz\label{surj}
\begin{lemm}
Consider the transposed map $N^t:T^*M\to T^*\F$ as a smooth map. Then for the canonical linear Poisson structure $\eta_{N}$ related to the Lie algebroid $T^*\F$ with the Lie algebra structure $[,]_{N}$ and the anchor $N$ the following equality holds:
$$
\eta_{N}=N^t_*\eta_{T^*M}.
$$
\end{lemm}
To prove this claim we shall proceed as in Remark \ref{remanc}.
First notice that for any $X\in\Ga(T\F)$ we have the following equality of linear functions on $T^*M$: $\overline{NX}=(N^t)^*\overline{X}$, where $(N^t)^*$ stands for the pullback. Now the following calculations, which use this equality and the definition of the algebroid $(T^*\F,[,]_{N},N)$, prove the claim (in view of properties (1)--(3) of the linear bracket which determine it, see Section \ref{sAl}):
\begin{align*}
(N^t)^*\{\overline X, \overline Y\}'=(N^t)^*\overline{[X,Y]_N}=\overline{N[X,Y]_N}=\overline{[NX,NY]}=\{\overline{NX},\overline{NY}\}
=\{(N^t)^*\overline{X},
(N^t)^*\overline{Y}\}\\
(N^t)^*\{\overline X,\sigma^* f\}'=(N^t)^*\sigma^*(NXf)=\pi^*(NXf)=\{\overline{NX},\pi^*f\}=\{(N^t)^*\overline{X},(N^t)^*\sigma^*f\}\\
(N^t)^*\{\sigma^*f,\sigma^* g\}'=0=\{\pi^*f,\pi^* g\}=\{(N^t)^*\sigma^*f,(N^t)^*\sigma^* g\};
\end{align*}
here $\{,\}'$ and $\{,\}$  are the Poisson brackets corresponding to $\eta_{N}$ and $\eta_{T^*M}$ correspondingly and  $\sigma:T^*\F\to M$ and $\pi:T^*M\to M$ are the canonical projections. \qed

\smallskip

Summarizing, the canonical  bihamiltonian structure on $T^*\F$ related to a PNO $(T\F,N)$ is generated by the linear Poisson structures $\eta_1:=(I^t)_*\eta_{T^*M}$ and $\eta_2:=(N^t)_*\eta_{T^*M}$. Note that the fibers of the canonical projection $T^*\F\to M$ are lagrangian submanifolds in any symplectic leaf of any of these two Poisson structures, i.e. the fibers form a {\em bilagrangian foliation}.

Below we consider two particular cases of the ``up construction''.

\abz\label{exaKro}
\begin{exa}\rm  Let $\{\F_s\}_{s\in\R\P^1}$ be a Kronecker web on a manifold
  $M$  and  $\phi_i\colon T^*M\to \Phi$ be the corresponding bundle morphisms (see Definition \ref{kronweb}).
In the particular case of the Kronecker PNO $(T\F_\infty,N)$, $N=\phi_2^t\circ(\phi_1^t)^{-1}$, related to a Kronecker web the ``up construction'' gives a  bihamiltonian structure $\eta_{1,2}\colon T^*M'\to TM'$, $M':=T^*\F_\infty$. We can say more about this bihamiltonian structure in comparison with the general case.

First of all, since $N$ is fiberwise injective, $N^t:T^*M\to T^*\F_\infty$ is a smooth surjective submersion and by Lemma \ref{surj} we can define $\eta_2=\eta_N$ as $N^t_*\eta_{T^*M}$.

 Second, let $x_1,\ldots,x_n$ be  a local coordinate system on $M$ such that $\frac{\d }{\d x_1},\ldots,\frac{\d }{\d x_k}$ are the basic vector fields tangent to $\F_\infty$ and let $\xi_1,\ldots,\xi_k$ be the corresponding linear functions on $T^*\F_\infty$. Then by formulas (\ref{linp}) the symplectic foliation $\overline{\F}_s$
 of the linear Poisson structure corresponding to the Lie algebroid $(T\F_\infty,N-s I)$ (here $I$ is the canonical embedding $T\F_\infty\hookrightarrow TM$) is generated by the vector fields $\frac{\d }{\d \xi_1},\ldots,\frac{\d }{\d \xi_k}$ and $(N-sI)\frac{\d }{\d \sigma^*x_1},\ldots,(N-sI)\frac{\d }{\d \sigma^*x_k}$ (here $\sigma: T^*\F_\infty\to M$ is the canonical projection, i.e $\sigma^*x_i$ is a base function on $T^*\F_\infty$). Due to the kroneckerity of $N$ the rank of the distribution $\overline D_s$ generated by these vector fields is constant even if we admit $s\in\C$, which means that the corresponding bihamiltonian structure $\eta_{1,2}$, is {\em Kronecker} itself, i.e.  for any
$p\in M'$ the J--K decomposition of the pair of operators $\eta_{1,p}, \eta_{2,p}\colon T_p^*M'\to T_pM'$ does not contain Jordan
blocks. Moreover, we observe the following obvious facts: (1) $\bigcap_{s}(\overline D_s)_p$ coincides with the fiber of $\sigma$ passing through $p\in M'$, i.e. the canonical bilagrangian foliation $\W_0$ of the Kronecker bihamiltonian structure $\eta_{1,2}$ (see Section \ref{constrA}) coincides with the foliation of fibers of $\sigma$; (2) the base of this foliation is correctly defined and coincides with $M$; (3) the projection of the symplectic foliation $\overline{\F}_s$ with respect to $\sigma$ coincides with the initial foliation $\F_s$ from the web for any $s$.
\end{exa}

\abz\label{Example 1'}
\begin{exa}\rm Let  $N\colon E\to TM$ be a PNO with the domain
  $E=TM$, i.e., $N$ is a ``usual" Nijenhuis operator. Then the
  ``up construction''  gives a bihamiltonian structure $\eta_1:=\eta_{T^*M}$, $\eta_2=\eta_N$ on the manifold
  $M':=T^*M$, where  $\eta_{T^*M}$ is the canonical Poisson structure on $M'=T^*M$.

  The (1,1)-tensor $N'\colon TM'\to TM'$ uniquely defined by $N'=\eta_N\circ\eta_{T^*M}^{-1}$ has zero Nijenhuis torsion due to the compatibility of $\eta_N$ and $\eta_{T^*M}$.   In case, when $N$ is fiberwise invertible, the Poisson structure $\eta_{N}$ is nondegenerate and $(N')^{-1}=\eta_{T^*M}\circ\eta_N^{-1}$ coincides with  the so-called cotangent lift of the operator $N$ defined as $\eta_{T^*M}\circ (N^t)^*\eta_{T^*M}^{-1}$ (see \cite{t3}); here the transposed operator $N^t\colon T^*M\to T^*M$ is
  regarded  as a smooth map of $M'$ and $\eta_{T^*M}^{-1}$ is the canonical symplectic form.  We know from Lemma \ref{surj} that the following equality holds $\eta_{N}:=N^t_*\eta_{T^*M}$, which  in the case of fiberwise invertible $N$ can serve as the definition of the linear Poisson structure $\eta_{N}$.

\end{exa}

\abz\label{typeN}
\begin{defi}\rm The bihamiltonian structure $\eta_{T^*M},\eta_N$ on $T^*M$ from Example \ref{Example 1'}  will be called the {\em bisymplectic} or {\em Jordan bihamiltonian structure of type $N$}.
\end{defi}
The last terminology is motivated by the fact that there are only Jordan blocks in the J--K decomposition of the pair of operators $\eta_{T^*M}|_p,\eta_N|_p:T^*_pM'\to T_pM'$ for any $p\in M'$.

\abz\label{remaB}
\begin{rema}\rm F. J. Turiel \cite{t3} proved that under some additional assumption of regularity (which is satisfied for generic cases) any Jordan bihamiltonian structure is locally equivalent to a bihmiltonian structure of type $N$. In the next section we shall also see that any Kronecker bihamiltonian structure is locally equivalent to the one built in Example \ref{exaKro}. Thus the examples above show that the notion of a PNO is a proper geometric framework for simultaneous treatment of Jordan and Kronecker bihamiltonian structures.
\end{rema}

\abz\label{realN}
\begin{lemm}
Let $(T\F,N)$ be a PNO on a manifold $M$. Assume there exists a Nijenhuis operator $\overline N:TM\to TM$ such that $N=\overline N|_{T\F}$.  Write $I:T\F\to TM$ for the canonical inclusion. Let $I^t:T^*M\to T^*\F$ be the transposed operator regarded as a smooth surjective submersion. Then $I^t_*\eta_{T^*M}=\eta_{T^*\F}$ and $I^t_*\eta_{\overline{N}}=\eta_{N}$.
\end{lemm}
The first equality was already discussed (see Remark \ref{remanc}). The second equality follows from the commutativity of the following diagram
$$
\xymatrix{
T\F \ar[r]^{I} \ar[dr]^{N} & TM \ar[d]^{\overline{N}}\\
& TM}
$$
(which implies $N^t= I^t\circ\overline{N}^t$ and in view of  Lemma \ref{surj} $\eta_{N}=(N^t)_*\eta_{T^*M}= I^t_*\circ(\overline{N}^t)_*\eta_{T^*M}=I^t_*\eta_{\overline{N}}$). \qed
%
\section{Relations of Kronecker webs with bihamiltonian structures}
\label{constrA}

There are two constructions relating Kronecker webs with bihamiltonian structures, which are mutually inverse in the sense that will be explained below (see \cite{gz2}, \cite{p1}, \cite{t2}).

Let $\eta_{1,2}\colon T^*M\to TM $ be a {\em Kronecker bihamiltonian structure}, i.e., a bihamiltonian structure such that for any
$x\in M$ the J--K decomposition of the pair of operators $\eta_{1,x}, \eta_{2,x}\colon T_x^*M\to T_xM$ does not contain Jordan
blocks.
The rank of the Poisson bivector $\la_1\eta_{1}+\la_2\eta_2$ does not depend on $x$ and $\la_{1,2}$ (when $(\la_1,\la_2)\neq0$);
denote by $\overline{\F}_\la$, $\la=\la_1:\la_2$, the corresponding symplectic foliation.
Then $\{\overline{\F}_\la\}_{\la\in\P^1}$ is a family of foliations of constant rank; as linear algebra shows, they contain a unique common subfoliation $\W_0$ such that $T_x\W_0=\bigcap_{\la\in\P^1}T_x\overline{\F}_\la$ for any $x\in M$. Such a foliation is lagrangian in any symplectic leaf of any of two Poisson structures and is called the \emph{bilagrangian foliation} of the Kronecker bihamiltonian structure.
Reduce attention to a sufficiently small open
subset $U \subset M$ on which the foliation $\W_0$ has a local base $\B$.

Finally, it turns out that $\B$ carries a rich geometric structure of a Kronecker web: a
collection of foliations in general position $\F_\la$ depending on $\la\in\P^1$ such that
the normal spaces $\N_m\F_\la\subset T^*_mM$ depend in a particular  way on parameter $\la$. These foliations are the ``projections'' of the foliations $\overline{\F}_\la$ w.r.t.~ the reduction of $U$ to $\B$. As we know from  Section \ref{VerKro}  such a structure is equivalent to  a geometric Kronecker PNO.

Note that the operators $\eta_{1,x}, \eta_{2,x}$ being skew symmetric necessarily contain both increasing and decreasing Kronecker blocks (see Definition \ref{10.10}) in the J--K decomposition, which are mutually transposed to each other. Algebraically the construction described, which will be referred to as \emph{``down construction''}, consists in cutting off the decreasing blocks.

Vice versa, let $\{\F_s\}_{s\in\R\P^1}$ be a Kronecker web on a manifold
  $M$  and  $\phi_i\colon T^*M\to \Phi$ be the corresponding bundle morphisms (see Definition \ref{kronweb}).
  Then ``up construction'', which was discussed in Section \ref{upCon}, applied to the Kronecker PNO $(T\F_\infty,N)$, $N=\phi_2^t\circ(\phi_1^t)^{-1}$, related to the Kronecker web  gives a Kronecker  bihamiltonian structure $\eta_{1,2}\colon T^*M'\to TM'$, $M':=T^*\F_\infty$.

From Example \ref{exaKro} we see that starting from a Kronecker web and applying first ``up construction'' and then ``down construction'' results in the initial Kronecker web.

Applying these constructions other way round is more subtle. Starting from any Kronecker bihamiltonian structure $\eta_{1,2}$ we can always perform locally ``down construction'' and get a Kronecker web $\{\F_s\}_{s\in\R\P^1}$. Applying to it the ``up construction'' results in a bihamiltonian structure $\eta'_{1,2}$ which a priori  need not coincide with the initial one. It was the initial conjecture of Gelfand and Zakharevich (formulated by them in the case of generic Kronecker bihamiltonian structures \cite{gz2}, i.e. with Kronecker webs which are Veronese webs) that the bihamiltonian structures $\eta_{1,2}$ and $\eta'_{1,2}$ are locally equivalent, i.e. there exists a local diffeomorfism bringing one structure to another.

This conjecture was proved by Turiel in the particular cases listed in the following theorem (see \cite[Theorem 3.2]{t8} and references therein).
\abz\label{turiel}
\begin{theo}(Turiel) A Kronecker bihamiltonian structure can be locally reconstructed from its Kronecker web obtined by means of the ``down construction'' in the following cases:
\begin{itemize}\item in complex or real analytic category;
\item in $C^\infty$ category for generic Kronecker bihamiltonian structures and Kronecker bihamiltonian structures with flat Kronecker webs.
\end{itemize}

\end{theo}
A Kronecker web $\{\F_s\}_{s\in\R\P^1}$ is called {\em flat} if in a vicinity of every point there exists a local diffeomorphism bringing simultaneously all the foliations $\F_s$ to the foliations of parallel planes on an open set in $\R^n$.

%
\section{Problem of local bisymplectic realization of a Kronecker bihamiltonian structure}
\label{realbi}

It is well known \cite{weinst} that, given a Poisson structure $\eta$ on a manifold $M$, for any point of $M$ there exists an open neighbourhood of this point $U$ and a symplectic manifold $(\overline U, \omega)$ with a surjective submersion $p:\overline U\to U$ such that $p_*\omega^{-1}=\eta|_U$; here $\omega^{-1}$ is the Poisson structure inverse to the symplectic form $\omega$. In other words, any Poisson structure has a  {\em local symplectic realization}. This is a first step to the problem of existence of {\em global} symplectic realization which  is very important and led in particular to the theory of symplectic groupoids.

Analogous problem can be formulated in the bihamiltonian context: given a bihamiltonian structure $\eta_{1,2}$ on a manifold $M$ such that $\la_1\eta_1+\la_2\eta_2$ is degenerate for any $\la$, does it have a bisymplectic realization, i.e. does there exist a manifold $\overline M$ with a bihamiltonian structure $\omega^{-1}_{1,2}$ (such bihamiltonian structures necessarily are  {\em Jordan}, i.e. for any $x\in M$ the pair of operators $\om^{-1}_{1,x}, \om^{-1}_{2,x}\colon T_x^*M\to T_xM$  contains only Jordan
blocks in the J--K decomposition) and a surjective submersion $p:\overline M\to M$ such that $p_*\omega^{-1}_{1,2}=\eta_{1,2}$? In this section we consider the problem of local bisymplectic realization for Kronecker bihamiltonian structures.

Note that there is a crucial difference between the two realization problems above: in the  Poisson case there is only one local model of the symplectic form $\om$ given by the Darboux theorem while there are many local models of bisymplectic bihamiltonian structures $\omega^{-1}_{1,2}$, i.e. Jordan bihamiltonian structures. For instance, the Jordan bihamiltonian structures of type $N$ (see Definition \ref{typeN}), which are completely determined by a Nijenhuis (1,1)-tensor $N$, are locally inequivalent for locally inequivalent  $N$.

A quite natural and desirable feature of the symplectic and bisymplectic realization is its minimality: once dimension of $M$ is fixed, try to find $\overline{M}$ of possibly minimal dimension. Since for a Kronecker bihamiltonian structure $\eta_{1,2}$ both the bivectors have the same rank, say $2r$, and corank, say $l$, it is easy to see that the minimal possible dimension for $\overline{M}$ we can think about is $2r+2l$.

Now we can make our problem more precise.

 \begin{probla}\begin{enumerate}\rm\item[(a)]
 Given a Kronecker bihamiltonian structure $\eta_{1,2}$, $\rk\eta_{1,2}=2r$, on an open set $U \subset M$, $\dim M=m$, do there exist a Jordan bihamiltonian structure $\overline{\eta}_{1,2}$ on an open set $\overline{U} \subset \overline{M}$, $\dim \overline{M}=2m-2r$, and a smooth surjective submersion $p:\overline{U}\to U$ such that $p_*\overline{\eta}_{1,2}=\eta_{1,2}$?
  \item[(b)]
  List all locally inequivalent Jordan bihamiltonian structures $\overline{\eta}_{1,2}$ on $\overline{U}$ with the property $p_*\overline{\eta}_{1,2}=\eta_{1,2}$.
  \end{enumerate}
\end{probla}

Below we set some preliminary steps for solving this problem. In view of Theorem \ref{turiel} we can assume that the bihamiltonian structure $\eta_{1,2}$ is equal to the bihamiltonian structure $\eta_{T^*\F_\infty},\eta_N$ on the manifold $T^*\F_\infty$, $T\F_\infty \subset T\B$, where $\B$ is the local base of the canonical bilagrangian foliation $\W_0$ of $\eta_{1,2}$ and  $N:T\F_\infty\to T\B$ is the Kronecker PNO corresponding to the Kronecker web obtained on $\B$ by means of the ``down construction'' (see Example \ref{exaKro}). Now assume that there exists a Nijenhuis operator $\overline N:T\B\to T\B$ such that $N=\overline N|_{T\F_\infty}$. Then by Lemma \ref{realN} we have
$$
I^t_*\eta_{T^*\B}=\eta_{T^*\F_\infty},I^t_*\eta_{\overline{N}}=\eta_{N},
$$
where $I:T\F_\infty \to T\B$ is the canonical inclusion and  $I^t:T^*\B\to T^*\F_\infty$ is the corresponding surjective submersion.

We see that Problem 1 is intimately related to the following

\begin{problb}
\begin{enumerate}\rm\item[(a)] Given a Kronecker PNO $(T\F,N)$, $\rk\F=r$, on an open set $V \subset \R^{m-r}$, $m>2r$, does there exist a Nijenhuis operator $\overline{N}:TV\to TV$ such that $\overline{N}|_{T\F}=N$?
\item[(b)] List all locally nonequivalent Nijenhuis operators $\overline{N}$ on $V$ satisfying $\overline{N}|_{T\F}=N$.
\end{enumerate}
\end{problb}

The considerations above show that once Problem 2(a) is solved we obtain also a  solution of Problem 1(a). Recall (see Remark \ref{realizA}) that Problem 2(a) has a solution for any Kronecker web, hence Problem 1(a) has a solution for any Kronecker bihamiltonian structure.

On the other hand, a solution of Problem 2(b), which will be called the \emph{``realization problem for  Kronecker webs''}, would imply only a particular solution of Problem 1(b), i.e. a solution in the class of Jordan bihamiltonian structures of type $\overline{N}$ on $\overline U$, where $\overline N:TV\to TV$ is a Nijenhuis (1,1)-tensor and $\overline U$ is an open set in $T^*V$ (cf. Remark \ref{remaB}). Solutions to the realization problem will be obtained in the next section for particular Kronecker webs.

%
\section{Realization problem for Veronese webs}
\label{realz}

The realization problem for Kronecker webs was formulated in the previous section (Problem 2(b)). Below we discuss this problem and we start from describing a solution to this problem for 3-dimensional Veronese webs obtained in \cite{kruglikovP}. We begin with the general situation, and then specify to the 3-dimensional case. For simplicity consider only \emph{complex analytic case} (which excludes  the normal form of a real Nijenhuis operator with complex eigenvalues, see \cite{kruglikovP} for this case).

Recall that one of the local models of the Nijenhuis operators $\overline{N}$, namely a semisimple operator with simple spectrum the elements of which are constant functions, was obtained in the proof of Theorem \ref{theo1q}. To get other local models we need to introduce the following notion.

 \abz\label{defSP}
\begin{defi}\rm
Consider a Veronese web $\{\F_\la\}_{\la\in\C\P^1}$ on a manifold $M^{n+1}$, given by $T\F_\la=\langle \al^\la\rangle^\bot$,
where $\al^\la=\al_0+\la\al_1+\cdots+\la^n\al_n$ and $\al_0,\al_1,\ldots,\al_n$ is a local coframe on an open set $U \subset M$.
An analytic  function $\phi:U\to \C$ is called {\em self-propelled}
if $d\phi$ is proportional to $\al^\phi$. If the coefficient of proportionality is nonzero, we denote this by $d\phi\sim\al^\phi$.
However, the coefficient is allowed to be zero, so a constant function is also considered self-propelled.
 \end{defi}

\abz\label{lemmSP}
\begin{lemm}
Let $\{\F_\la\}$ be a Veronese web on $M^{n+1}$. Then in a vicinity of any point $x\in M$ there exist $n+1$ functionally independent self-propelled
functions $\phi_0(x),\phi_1(x),\dots,\phi_n(x)$. If $X_0,\dots,X_n$ is the frame dual to the coframe $\al_0,\dots,\al_n$ defining the Veronese web,
the condition on the function $\phi$ to be self-propelled is the following system PDEs:
 \begin{equation}\equ\label{sys}
\phi X_0\phi=X_1\phi,\ldots,\phi X_{n-1}\phi=X_n\phi.
 \end{equation}
\end{lemm}

\noindent The required relation $\al_0+\dots+\phi^n\al_n\sim (X_0\phi) \al_0+\dots+(X_n\phi) \al_n$ is equivalent to
vanishing of the determinants
$$
\left|
  \begin{array}{cc}
  1 & \phi \\ X_0\phi  & X_1\phi  \\
  \end{array}
\right|,\dots,\
\left|
  \begin{array}{cc}
  \phi^{n-1} & \phi^n \\ X_{n-1}\phi  & X_n\phi  \\
  \end{array}
\right|,
$$
which is equivalent to system (\ref{sys}). Let $F(x,\la)$ be a $\la$-parametric first integral of the folitation $\F_\la$, where  $x=(x_1,\dots,x_n)$.
The following formula gives a family of implicit solutions $\phi(x)$ of system \eqref{sys} depending on an arbitrary smooth function of one
variable $f=f(\la)$ that locally satisfies $f'(\la)\not=F_\la$:
 \begin{equation}\equ\label{eqt}
F(x,\phi(x))=f(\phi(x)).
 \end{equation}
Indeed, differentiating this equality along $X_k-\phi(x) X_{k-1}$  we get
 \begin{equation}\equ\label{eqt1}
d_xF(x,\lambda)(X_k-\la X_{k-1})|_{\lambda=\phi(x)}+(F_\la(x,\phi(x))-f'(\phi(x)))\cdot(X_k\phi(x)-\phi(x) X_{k-1}\phi(x))=0.
 \end{equation}
The first term vanishes since $X_k-\la X_{k-1}\in\langle\al^\la\rangle^\perp$, and the claim follows.

Choosing $n$ solutions $\phi_0,\ \dots,\ \phi_n$ with initial values $c_0,\dots,c_n$ at $x\in M$ being
pairwise different and with nonzero $\psi_i:=X_0\phi_i|_x$, we compute from (\ref{sys}) the Jacobian at $x$:
 $$
\mathrm{Jac}_{x}(\phi_0,\phi_1,\dots,\phi_n)\sim
\left|\begin{array}{cccc}\psi_0 & c_0\psi_0 & \dots & c_0^n\psi_0\\ \psi_1 & c_1\psi_1 & \dots & c_1^n\psi_1\\
\vdots & \vdots & \ddots & \vdots \\
\psi_n & c_n\psi_n & \dots & c_n^n\psi_n\end{array}\right|
=\psi_0\psi_1\cdots\psi_n\left|
\begin{array}{cccc}1 & c_0 & \dots & c_0^n\\ 1 & c_1 & \dots & c_1^n\\
\vdots & \vdots & \ddots & \vdots\\ 1 & c_n & \dots & c_n^n\end{array}\right|.
 $$
Since the Vandermonde determinant with the second column consisting of pairwise different entries is nonzero,
we obtain $n$ functionally independent solutions of (\ref{sys}). \qed

\abz\label{mainTh}
\begin{theo}
Let $(T\F,N)$ be a Kronecker PNO of generic type (see Theorem \ref{theo1q}) on a 3-dimensional manifold $M$. Then in a
neighborhood $U$ of every point $p\in M$ there exists a Nijenhuis operator $\overline{N}:TM\to TM$ of any type A, B or C listed in Appendix  such that $\overline{N}|_{T\F}=N$.
\end{theo}

\noindent Consider $(T\F, N)$ locally near $p\in M$. The intersection $D_1:=T\F\cap NT\F$ is a one dimensional distribution. Choose arbitrarily a nonvanishing vector field $X_1\in \G(D_1)$ and put $X_0:=N^{-1}X_1,X_2:=NX_1$. Then $X_0,X_1,X_2$ is a frame such that there exist functions $b_0,b_1,c_1,c_2$ satisfying the following commutation relations:
\begin{itemize}
\item[(i)]  $[X_0,X_1]=b_0X_0+b_1X_1$ and $ [X_1,X_2]=c_1X_1+c_2X_2$;
    \item[(ii)] $[X_0,X_2]=c_1X_0+(c_2+b_0)X_1+b_1X_2$.
     \end{itemize}
Item (i) is  due to the integrability of the distributions $T\F$ and $NT\F$.    To prove Item (ii) let $[X_0,X_1]=d_0X_0+d_1X_1+d_2X_2$ for some functions $d_0,d_1$, and $d_2$ and use
the definition of a PNO \ref{PNOdefi}: by condition 1 of this definition we have $[X_0,X_1]_N=$ $[NX_0,X_1]+$ $[X_0,NX_1]-$ $N[X_0,X_1]=[X_1,X_1]+[X_0,X_2]-N[X_0,X_1]=[X_0,X_2]-N[X_0,X_1]=d_0X_0+d_1X_1+d_2X_2-
(b_0X_1+b_1X_2)=d_0X_0+(d_1-b_0)X_1+(d_2-b_1)X_2\in T\F$, which implies $d_2=b_1$; by condition 2 of this definition we have
$c_1X_1+c_2X_2=[X_1,X_2]=[NX_0,NX_1]=N([X_0,X_1]_N)=N(d_0X_0+(d_1-b_0)X_1)=
d_0X_1+(d_1-b_0)X_2$, which implies $d_0=c_1, d_1=c_2+b_0$.

If $(X_0,X_1,X_2)$ is a frame satisfying relations (i-ii) for some functions and $(\al_0,\al_1,\al_2)$ is the dual coframe, it is easy to see that the distribution $\langle \al_0+\la\al_1+\la^2\al_2\rangle^\bot \subset TM$ is integrable for any $\la$, i.e. defines a Veronese web $\{\F_\la\}$. This is of course the Veronese web corresponding to $N$ by Theorem \ref{theo1q} (see its proof).

The matrix of the operator $N:T\F\to TM$ with respect to the bases $(X_0,X_1)$ in $T\F$ and $(X_0,X_1,X_2)$ in $TM$ is equal to
$$
\left[
  \begin{array}{cc}
    0 &   0  \\
    1 & 0  \\
    0 & 1  \\
  \end{array}
\right].
$$
Define $\overline{N}$ by $\overline{N}|_{T\F}=N$ and $\overline{N}X_2=f_0 X_0+f_1X_1+f_2X_2$, where $f_i$ are local analytic functions,
i.e. putting the matrix of $\overline{N}$ in the frame $(X_0,X_1,X_2)$ to be equal to
$$
\left[
  \begin{array}{ccc}
    0 & 0 & f_0 \\
    1 & 0 & f_1 \\
    0 & 1 & f_2 \\
  \end{array}
\right].
$$
Direct calculations taking into account relations (i), (ii) show that
$T_{\overline{N}}(X_1,X_2)=0$, if and only if the following system of nonlinear first order equations is satisfied:
\begin{equation}\equ\label{er1}
X_2f_0=f_0X_1f_2, X_2f_1=X_1f_0+f_1X_1f_2, X_2f_2=X_1f_1+f_2X_1f_2,
\end{equation}
and, analogously, the equality
$T_{\overline{N}}(X_0,X_2)=0$
is equivalent to the system
\begin{equation}\equ\label{er2}
X_1f_0=f_0X_0f_2, X_1f_1=X_0f_0+f_1X_0f_2, X_1f_2=X_0f_1+f_2X_0f_2.
\end{equation}
Now let $f_1=\phi_1\phi_2\phi_3, f_2=-\phi_1\phi_2-\phi_1\phi_3-\phi_2\phi_3, f_3=\phi_1+\phi_2+\phi_3$ for some  local functions $\phi_1,\phi_2,\phi_3$. Then it is easy to see that once the functions $\phi_i$ satisfy the system of equations (\ref{sys}), the functions $f_i$ satisfy  the systems of equations (\ref{er1}), \ref{er2}). In other words,  if the functions $\phi_1,\phi_2,\phi_3$ are self-propelled for the corresponding Veronese web, the Nijenhuis torsion $T_{\overline{N}}$ of the (1,1)-tensor $\overline{N}$ given in the frame $X_0,X_1,X_2$ by the matrix
\begin{equation}\equ\label{matr}
F(\phi_1,\phi_2,\phi_3):=\left[
  \begin{array}{ccc}
  0 & 0 &  \phi_1\phi_2\phi_3\\
    1 & 0 & -\phi_1\phi_2-\phi_1\phi_3-\phi_2\phi_3 \\
    0 & 1 & \phi_1+\phi_2+\phi_3 \\
  \end{array}
\right]
\end{equation}
vanishes (recall that $T_{\overline{N}}(X_0,X_1)=T_{N}(X_0,X_1)=0$ by the assumptions of the theorem).

Now let $\psi_1,\psi_2,\psi_3$ be functionally independent self-propelled functions with pairwise distinct $\psi_1(p)$, $\psi_2(p)$, $\psi_3(p)$ (they exist by Lemma \ref{lemmSP}) and let $a_1,a_2,a_3$ be pairwise distinct constants. Put
\begin{itemize}
\item $F_{A0}:=F(\psi_1,\psi_2,\psi_3)$; $F_{A1}:=F(\psi_1,\psi_2,a_3)$; $F_{A2}:=F(\psi_1,a_2,a_3)$; $F_{A3}:=F(a_1,a_2,a_3)$;
\item  $F_{B0}:=F(\psi_2,\psi_2,\psi_3)$; $F_{B1}:=F(\psi_2,\psi_2,a_3)$; $F_{B2}:=F(a_2,a_2,\psi_3)$; $F_{B3}:=F(a_2,a_2,a_3)$;
\item  $F_{C0}:=F(\psi_3,\psi_3,\psi_3)$; $F_{C1}:=F(a_3,a_3,a_3)$.
\end{itemize}

We have shown above that all these matrices represent Nijenhuis (1,1)-tensors.  On the other hand, we recognize in these matrices the Frobenius forms of all the Nijenhuis (1,1)-tensors listed in Appendix. Consequently, by \cite{t4} for each $F_X$ there should exist local coordinates $(x_1,x_2,x_3)$ such that the matrix of  the corresponding Nijenhuis (1,1)-tensor $\overline{N}$ in the basis $\{\frac{\d }{\d x_i}\}$ has the form $N_X$ from the list of Appendix. \qed

\smallskip

We conclude this section by a conjecture that the realization problem for a Veronese web can be similarly solved in any dimension (in fact its proof should go in the same way as above).

\abz\label{conj}
\begin{conj} Let $(T\F,N)$ be a Kronecker PNO of generic type (see Theorem \ref{theo1q}) on a $n$-dimensional manifold $M$, $n>3$. Then in a
neighborhood $U$ of every point $p\in M$ there exists a cyclic Nijenhuis operator $\overline{N}:TM\to TM$ of any type of \cite{t4} such that $\overline{N}|_{T\F}=N$.
\end{conj}
Note that condition of cyclicity is necessary when we speak about the extensions of PNOs of generic type.

\abz\label{remn}
\begin{rema}\rm
The situation with nongeneric Kronecker PNOs, i.e. having more than one Kronecker block in the J--K decomposition seems to be much more involved. The extension here can have more than one cyclic blocks, however not necessarily. The analogues of systems of equations (\ref{sys}), (\ref{er1}), (\ref{er2}) would be much more complicated.
\end{rema}

%
\section{The Hirota equation}
\label{Hiroo}

In this section we assume that $\dim M=3$.
The aim of this section is to show that there is  a 1--1-correspondence between Veronese webs in 3  dimensions and solutions of the so-called dispersionless Hirota PDE
$$
a_1f_{x_1}f_{x_2x_3}+a_2f_{x_2}f_{x_3x_1}+a_3f_{x_3}f_{x_1x_2}=0,
$$
where $a_i$ are constants such that $a_1+a_2+a_3=0$.

It follows from Theorem \ref{theo1q} and its proof that, given a Veronese web, one can construct a PNO which, at least locally, can be extended to a
Nijenhuis operator defined on the whole tangent bundle $TM$.
 In Section \ref{realz} we have shown that in fact such an extension is possible essentially to any of normal forms of Nijenhuis operators in 3 dimensions.

Conversely, starting from a Nijenhuis (1,1)-tensor ${N}$  we can try to construct a Veronese web by means of constructing a PNO (cf. Theorem \ref{theo1q}) $(\F,{N}|_{T\F})$  for some foliation $\F$. Assuming that the foliation $\F$ is given by $f=\op{const}$ for some smooth function $f$, we can use Lemma \ref{llll} to obtain sufficient conditions for ${N}|_{T\F}$
to be a PNO in terms of a PDE on $f$, the form of which essentially depends on the form of the initial Nijenhuis operator.

Let us illustrate these idea choosing the simplest normal form of a Nijenhuis operator: the diagonal one with constant pairwise distinct eigenvalues.

\abz\label{constr}
\begin{constr}\rm
Consider $M=\R^3(x_1,x_2,x_3)$ and
a Nijenhuis operator $N:TM\to TM$ defined by
 \begin{equation}\equ\label{Jx}
N\partial_{x_i}=\la_i\partial_{x_i},
 \end{equation}
 where $\la_1,\la_2,\la_3$ are pairwise distinct nonzero numbers.
Let $f:\R^3\to \R$ be a smooth function such that $f_{x_i}\not=0$. Define a foliation $\F_\infty$ by $f=\op{const}$, i.e. by
$T\F_\infty:=\langle df\rangle^\bot$. Then $(N(T\F_\infty))^\bot=\langle\om\rangle$, where
$$
\om=(N^t)^{-1}df= \la_1^{-1}f_{x_1}dx_1+\la_2^{-1}f_{x_2}dx_2+\la_3^{-1}f_{x_3}dx_3.
$$
The condition of integrability of the distribution $N(T\F_\infty)$, $d\om\wedge\om=0$ (which by Lemma \ref{llll} implies that $N|_{T\F_\infty}$ is a PNO), is equivalent to
 \begin{equation}\equ\label{Hirota-a}
(\la_2-\la_3)f_{x_1}f_{x_2x_3}+(\la_3-\la_1)f_{x_2}f_{x_3x_1}+(\la_1-\la_2)f_{x_3}f_{x_1x_2}=0,
 \end{equation}
in which we recognize the Hirota equation.
\end{constr}
The following theorem is a variant of \cite[Theorem 3.8]{z2} (our proof is different).

\abz\label{6}
\begin{theo} Let $\la_1,\la_2,\la_3$ be  distinct real numbers.
 \begin{enumerate}
\item For any solution $f$ of (\ref{Hirota-a}) on a domain $U\subset M$ with $f_{x_i}\not=0$, $i=1,2,3$, the 1-form
 \begin{equation}\equ\label{1-f}
\al^\la=(\la_2-\la)(\la_3-\la)f_{x_1}dx_1+(\la_3-\la)(\la_1-\la)f_{x_2}dx_2+(\la_1-\la)(\la_2-\la)f_{x_3}dx_3
 \end{equation}
defines a Veronese web $\F_\la$ on $U$ by $T\F_\la=\langle\al^\la\rangle^\bot$ such that
 \begin{equation}\equ\label{FFF}
\F_{\la_i}=\{x_i=const\},\ \F_\infty=\{f=const\}.
 \end{equation}
\item Conversely, let $\{\F_\la\}$ be a Veronese web on a 3-dimensional smooth manifold $M$. Then in a neighbourhood of any point on $M$ there exist local coordinates $(x_1,x_2,x_3)$ such that any smooth first integral $f$ of the foliation $\F_\infty$ is a solution of equation (\ref{Hirota-a}) with $f_{x_i}\not=0$.
 \end{enumerate}
Consequently, we obtain a 1--1-correspondence between Veronese webs $\{\F_\la\}$ satisfying (\ref{FFF}) and
the classes $[f]$ of solutions $f$ of (\ref{Hirota-a}) with  $f_{x_i}\not=0$ modulo the following equivalence relation:
$f\sim g$ if there exist local diffeomorphisms $\psi,\phi_1,\phi_2,\phi_3$ of $\R$ such that $f(x_1,x_2,x_3)=\psi(g(\phi_1(x_1),\phi_2(x_2),\phi_3(x_3))$
(obviously, if $f\sim g$ and $f$ solves (\ref{Hirota-a}), then $g$ does the same).
 \end{theo}

\noindent
On a solution $f$ of equation (\ref{Hirota-a}) we get $d\om\wedge\om=0$, hence the distribution $N(T\F_\infty)$ is integrable.
Consequently, $N|_{T\F_\infty}$ is a PNO by Lemma \ref{llll}. The condition $f_{x_i}\not=0$
implies that the pair $(N|_{T\F_\infty},I)$ has generic type (one Kronecker block in the J--K decomposition) and thus defines a Veronese web $\F_\la$ by Theorem \ref{theo1q}.
The Veronese curve $\al^\la$ in $T^*U$ such that $(T\F_\la)^\bot=\langle \al^\la \rangle$ annihilates the distribution $N_\la(T\F_\infty)=T\F_\la$.
Direct check shows that it is given by formula (\ref{1-f}), in particular satisfies (\ref{FFF}).

Conversely, let $\F_\la$ be a Veronese web and $f$ a first integral of $\F_\infty$. The proof of Theorem \ref{theo1q} yields the coordinates
$(x_1,x_2,x_3)$ and a Nijenhuis operator by (\ref{Jx}). The distribution $N(T\F_\infty)=T\F_0$ is integrable, hence $d\om\wedge\om=0$ and $f$
solves (\ref{Hirota-a}). The condition $f_{x_i}\not=0$ follows from nondegeneracy of the curve $\al^\la$.

Finally, the last statement follows from the fact that the first integrals of the three Veronese foliations corresponding to different $\la_1,\la_2,\la_3$
determine the first integral of any other foliation up to postcomposition with a local diffeomorphism.
 \qed

\bigskip

%
\section{Other PDEs related to Veronese webs in 3D}
\label{listPDE}

Repeating Construction \ref{constr} for other types of Nijenhuis operators listed in Appendix we get another PDEs on the function $f$, which are pairwise contactly nonequivalent (see \cite[Section 6]{kruglikovP}). Below we list these PDEs corresponding to the cases A, B, C of Appendix, and indicate the Veronese curves $\al^\la$
(the the one-forms $\omega$ such that $(N(T\F_\infty))^\bot=\langle\om\rangle$ are given   by $\om=\al^\la|_{\la=0}$).
\begin{enumerate}
\item[\bf(A)]\hspace{0.5cm}
\framebox{ $(\la_2(x_2)-\la_3(x_3))f_{x_1}f_{x_2x_3}+(\la_3(x_3)-\la_1(x_1))f_{x_2}f_{x_3x_1}+(\la_1(x_1)-\la_2(x_2))f_{x_3}f_{x_1x_2}=0$}
 $$
\hspace{-0.7cm}\al^\la = (\la_2(x_2)-\la)(\la_3(x_3)-\la)f_{x_1}dx_1+
(\la_3(x_3)-\la)(\la_1(x_1)-\la)f_{x_2}dx_2+(\la_1(x_1)-\la)(\la_2(x_2)-\la)f_{x_3}dx_3.
 $$
\item[\bf(B)]\hspace{0.5cm}
\framebox{ $f_{x_1}f_{x_1x_3}-f_{x_3}f_{x_1x_1}+(\la_2(x_2)-\la_3(x_3))(f_{x_1}f_{x_2x_3}-f_{x_2}f_{x_1x_3})
+\la_2'(x_2)f_{x_1}f_{x_3}=0$}
 $$
\hspace{-0.7cm}\al^\la = (\la_2(x_2)-\la)(\la_3(x_3)-\la)(f_{x_1}dx_1+f_{x_2}dx_2)+(\la_2(x_2)-\la)^2f_{x_3}dx_3
-(\la_3(x_3)-\la)f_{x_1}dx_2.
 $$
\item[\bf (C)] C0 \hspace{1cm}  \framebox{ $(f_{x_1}f_{x_2x_2}-f_{x_2}f_{x_1x_2})x_2+f_{x_3}f_{x_2x_2}-f_{x_2}f_{x_2x_3}+
f_{x_2}f_{x_1x_1}-f_{x_1}f_{x_1x_2}+f_{x_1}f_{x_2}=0$}
\begin{eqnarray*}
\al^\la = f_{x_1}((x_3-\la)^2dx_1-(x_3-\la)dx_3)+
f_{x_2}(-(x_3-\la)dx_1+(x_3-\la)^2dx_2\\
+(x_2(x_3-\la)+1)dx_3)+  f_{x_3}(x_3-\la)^2dx_3.
\end{eqnarray*}
\item[] C1\hspace{1cm}  \framebox{ $f_{x_1}f_{x_3x_1}-f_{x_3}f_{x_1x_1}+f_{x_2}f_{x_1x_2}-f_{x_1}f_{x_2x_2}=0$}
\begin{eqnarray*}
\al^\la = f_{x_1}((a_3-\la)^2dx_1-(a_3-\la)dx_2+dx_3)+ f_{x_2}((a_3-\la)^2dx_2\\
-(a_3-\la)dx_3)+f_{x_3}(a_3-\la)^2dx_3.
\end{eqnarray*}
 \end{enumerate}
Here the following specifications should be made in order to exhaust the corresponding cases  ($a_1,a_2,a_3$ are arbitrary pairwise different constants):
 \begin{enumerate}
\item[\bf(A)]
$A0$: $\la_1(x_1)=x_1, \la_2(x_2)=x_2, \la_3(x_3)=x_3$;\qquad
$A1$: $\la_1(x_1)=x_1, \la_2(x_2)=x_2, \la_3(x_3)=a_3$;\\
$A2$: $\la_1(x_1)=x_1, \la_2(x_2)=a_2, \la_3(x_3)=a_3$;\qquad\
$A3$: $\la_1(x_1)=a_1, \la_2(x_2)=a_2, \la_3(x_3)=a_3$.
\item[\bf(B)]
$B0$: $\la_2(x_2)=x_2, \la_3(x_3)=x_3$;\qquad
$B1$: $\la_2(x_2)=x_2, \la_3(x_3)=a_3$;\\
$B2$: $\la_2(x_2)=a_2, \la_3(x_3)=x_3$;\qquad\,%
$B3$: $\la_2(x_2)=a_2, \la_3(x_3)=a_3$.
 \end{enumerate}

Note that case $A3$ corresponds to the Hirota equation considered in the previous section. It turns out that in fact for all these equations there is a 1--1-correspondence between Veronsese webs and classes of ``nondegenerate'' solutions with respect to the natural equivalence relation, i.e. the analogue of Theorem \ref{6} holds; here the solutions are nondegenerate in the following sense.
\abz\label{defiND}
\begin{defi}\rm
A solution $f$ of any of the equations A, B, C on an open set $U\subset M$ with coordinates $(x_1,x_2,x_3)$ is called {\em nondegenerate} if the
corresponding one-form $\al^\la\in T^*U$ defines a Veronese curve at any $x\in U$ (equivalently: the curve
$\la\mapsto\al^\la=\al_0+\la\al_1+\la^2\al_2$ does not lie in any plane, i.e., the 1-forms $\al_0, \al_1, \al_2$ are linearly independent at any point).
\end{defi}

\abz\label{theoP1}
\begin{theo}
\begin{enumerate}
\item A generic solution $f$ of any of the equations A, B, C is nondegenerate on a small open set $U$. If $f$ is such a solution,
then the corresponding one-form $\al^\la$ defines a Veronese web $\F_\la$ on $U$ by $T\F_\la=\langle\al^\la\rangle^\bot$.
\item Conversely, let $\F_\la$ be a Veronese web on a 3-dimensional smooth manifold $M$. Then for any symbol $S=Ai,Bi,Ci$ in a neighbourhood of any point on $M$ there exist local coordinates $(x_1,x_2,x_3)$ such that any smooth first integral $f$ of the foliation $\F_\infty$ is a nondegenerate solution of the equation of type $S$.
\end{enumerate}
\noindent Here by a generic solution we mean a solution with a generic jet in the Cauchy problem setup. We omit the formulation of the analogue of the last part of Theorem \ref{6} as it follows immediately.
\end{theo}

\medskip

\noindent
The proof of the second statement of Item 1 is the same as that of Theorem \ref{6}(1).
For the explanation why a generic solution of the equations A--C is nondegenerate see \cite[Theorem 5.2]{kruglikovP}.

The  proof Item 2 goes essentially as that of Theorem \ref{6}(2) with the account of Theorem \ref{mainTh}.
\qed

\section{Generalizations to higher dimensions: systems of PDEs}
\label{high}

Generalization of the correspondence between Veronese webs and PDEs to higher dimensions (and the case of Kronecker webs) is straightforward. Let $\{\F_\la\}$ be a Kronecker web defined on an open set $U \subset \R^n$ and let  $N:T\F_\infty\to TU$ be the corresponding kronecker PNO, see Remark \ref{kronPNO}. By Remark \ref{realizA} (see also Theorem \ref{mainTh}) there exists a Nijehuis operator $\overline{N}:TU\to TU$   such that $\overline{N}|_{T\F_\infty}=N$. If $f_1,\ldots,f_k$ are functionally independent first integrals of the foliation $\F_\infty$, the condition of the integrability of the distribution $\overline N(T\F_\infty)$ (which follows from Lemma \ref{GPNO.10}) is equivalent to a system of nonlinear PDEs on the functions $f_i$ (depending on the form of the extension $\overline N$).

Conversely, given a Nijenhuis (1,1)-tensor $\overline N:TU\to TU$ and foliation $\F$ on $U$ defined by a system of first integrals $f_1,\ldots,f_k$, we can try to construct a Kronecker PNO $(T\F,\overline N|_{T\F})$ (and thus a Kronecker web $\{\F_\la\}$ with $\F_\infty=\F$ by requiring the integrability of the distribution $\overline N(T\F)$ (cf. Lemma \ref{llll}). The condition of the integrability of $\overline N(T\F)$ is equivalent to  a system of nonlinear PDEs on the functions $f_i$. Of course, one should impose additional algebraic conditions on the pair $\F,\overline N$ in order to guarantee the kroneckerity of $N|_{T\F}$.

Note that  the system of PDEs mentioned is overdetermined unless rank of the foliation $\F$ is not equal two (a reasonable bound is $\rk\F\ge 2$, since the rank one case gives a trivial differential constraint). Say in the case of Veronese web in $4$ dimensions we get $4$ equations on one function (the components of the $3$-form $d\om\wedge\om$, where $\om$ is the $1$-form annihilating $\overline N(T\F)$).

Let us illustrate the simplest higher dimensional case, when the system is determined: a Kronecker web $\{\F_\la\}$ with foliations $\F_\la$ of rank two in $M=\R^4$. If $N:T\F_\infty\to TM$ is the corresponding PNO, the pair $(N,I)$ has two Kronecker blocks in the J--K decomposition. It is known that such Kronecker webs are related with  torsionless $3$-webs on $M$, i.e. triples of foliations of rank $2$ in general position with the torsionless Chern connection. For any $3$-web $(\F^1,\F^2,\F^3)$ there exists a unique $1$-parametric family of distributions $\{D_\la\}_{\la\in\R\P^1}$ of rank two such that $D_\infty=T\F^1$, $D_0=T\F^2$, $D_1:=T\F^3$ and $D_\la$ is integrable for any $\la$ if and only if the torsion of the canonical Chern connection vanishes \cite[Theorem 4.14]{Nagy}. The corresponding family of foliations $\{\F_\la\}$ form a Kronecker web. We shall say that such Kronecker webs are {\em of $3$-web type}.

Consider $M=\R^4(x_1,x_2,x_3,x_4)$ and
a Nijenhuis operator $N:TM\to TM$ defined by
 \begin{equation}
N\partial_{x_i}=\la_i\partial_{x_i},\nonumber
 \end{equation}
 where $\la_1,\la_2,\la_3,\la_4$ are pairwise distinct nonzero numbers.
Let $f^{1,2}:\R^4\to \R$ be a generic pair of smooth functions. Define a foliation $\F_\infty$ by
$T\F_\infty:=\langle df^1, df^2\rangle^\bot$. Then $(N(T\F_\infty))^\bot=\langle\om_1,\om_2\rangle$, where
$$
\om_i=(N^t)^{-1}df^i= \sum_{j=1}^4\la_j^{-1}f^i_{x_j}dx_j.
$$
The condition of integrability of the distribution $N(T\F_\infty)$,
$$
d\om_1\wedge\om_1\wedge\om_2=0,d\om_2\wedge\om_1\wedge\om_2=0,
$$
 is equivalent to the following system of equations
 \begin{eqnarray}\equ\label{Hirot} \nonumber
(\la_1-\la_2)f^i_{x_1x_2}(f^i_{x_3}f^{\bar{i}}_{x_4}-f^{\bar{i}}_{x_3}f^i_{x_4})+
(\la_3-\la_1)f^i_{x_1x_3}(f^i_{x_2}f^{\bar{i}}_{x_4}-f^{\bar{i}}_{x_2}f^i_{x_4})&+&\\ \nonumber
(\la_1-\la_4)f^i_{x_1x_4}(f^i_{x_2}f^{\bar{i}}_{x_3}-f^{\bar{i}}_{x_2}f^i_{x_3})+
(\la_2-\la_3)f^i_{x_2x_3}(f^i_{x_1}f^{\bar{i}}_{x_4}-f^{\bar{i}}_{x_1}f^i_{x_4})&+&\\
(\la_4-\la_2)f^i_{x_2x_4}(f^i_{x_1}f^{\bar{i}}_{x_3}-f^{\bar{i}}_{x_1}f^i_{x_3})+
(\la_3-\la_4)f^i_{x_3x_4}(f^i_{x_1}f^{\bar{i}}_{x_2}-f^{\bar{i}}_{x_1}f^i_{x_2})&=&0,\ i=1,2,
 \end{eqnarray}
where $\bar{1}=2$ and $\bar{2}=1$. Once this system is satisfied by a pair of functions $f^{1,2}$, the distribution $D_\la=(N-\la\Id_{TM})(T\F_\infty)$  is integrable for any $\la$ and generates a Kronecker web, as it follows from Lemma \ref{llll}. Note that $D_\la$ is
annihilated by the pair of forms $\om_{1,2}^\la:=\sum_{j=1}^4(\la-\la_1)\cdots\widehat{(\la-\la_j)}\cdots(\la-\la_4)f^{1,2}_{x_j}dx_j$ for any $\la\not=\la_j$, where $\widehat{(\cdot)}$ means that the corresponding term is omitted (for $\la=\la_j$ the 1-forms $\om_{1,2}^\la$ become linearly dependent).

Another system of equations is obtained if we consider a Nijenhuis operator with two double eigenvalues. For instance, put $\la_1=\la_2$ and $\la_3=\la_4$ in the example above. Then  system (\ref{Hirot}) becomes
 \begin{eqnarray}\equ\label{Hiro}
f^i_{x_1x_3}(f^i_{x_2}f^{\bar{i}}_{x_4}-f^{\bar{i}}_{x_2}f^i_{x_4})-
f^i_{x_1x_4}(f^i_{x_2}f^{\bar{i}}_{x_3}-f^{\bar{i}}_{x_2}f^i_{x_3})&-&\\ \nonumber
f^i_{x_2x_3}(f^i_{x_1}f^{\bar{i}}_{x_4}-f^{\bar{i}}_{x_1}f^i_{x_4})+
f^i_{x_2x_4}(f^i_{x_1}f^{\bar{i}}_{x_3}-f^{\bar{i}}_{x_1}f^i_{x_3})
&=&0,\ i=1,2,
 \end{eqnarray}
and the corresponding annihilating one-forms are
\begin{equation}\equ\label{erp}
\om^\la_{1,2}=(\la-\la_1)(f^{1,2}_{x_1}dx_1+f^{1,2}_{x_2}dx_2)+
(\la-\la_3)(f^{1,2}_{x_3}dx_3+f^{1,2}_{x_4}dx_4)
 \end{equation}
 (now they span $D_\la^\bot$ for all $\la$). We see that $D_{\la_i}$, $i=1,3$, coincide with the corresponding coordinate planes. We can prove an analogue of Theorem  \ref{6}.

\abz\label{th3web}
\begin{theo}
Let $\la_1,\la_3$ be  distinct real numbers.
 \begin{enumerate}
\item For any solution $f^{1,2}$ of (\ref{Hiro}) on a domain $U\subset M$ satisfying
\begin{equation}\equ\label{jac}
 \left|\frac{D(f^1,f^2)}{D(x_1,x_2)}\right|\not=0,\left|\frac{D(f^1,f^2)}{D(x_3,x_4)}\right|\not=0
  \end{equation}
  the 1-forms (\ref{erp})
define a Kronecker web  $\F_\la$ on $U$ of $3$-web type by $T\F_\la=\langle\om^\la_1,\om^\la_2\rangle^\bot$ such that
 \begin{equation}\equ\label{FFFq}
\F_{\la_i}=\{x_i=const,x_{i+1}=const\},\ \F_\infty=\{f^1=const,f^2=const\}.
 \end{equation}
 \item
Conversely, let $\{\F_\la\}$ be a
 Kronecker  web of $3$-web type on a 4-dimensional smooth manifold $M$. Then in a neighbourhood of any point on $M$ there exist local coordinates $(x_1,x_2,x_3,x_4)$ such that any independent smooth first integrals $f^{1,2}$ of the foliation $\F_\infty$ are solutions of system (\ref{Hiro}) satisfying (\ref{jac}).
 \end{enumerate}
Consequently, we obtain a 1--1-correspondence between Kronecker  webs  $\{\F_\la\}$ of $3$-web type satisfying (\ref{FFFq}) and
the classes $[f^{1,2}]$ of solutions $f^{1,2}$ of (\ref{Hiro}) satisfying system (\ref{jac}) modulo the following equivalence relation:
$f^{1,2}\sim g^{1,2}$ if there exist local diffeomorphisms $\psi=(\psi_1,\psi_2),\phi=(\phi_1,\phi_2)$ and $\zeta=(\zeta_1,\zeta_2)$ of $\R^2$ such that
\begin{eqnarray*}
f^{1,2}(x_1,x_2,x_3,x_4)=\\
\psi_{1,2}(g^1(\phi_1(x_1,x_2),\phi_2(x_1,x_2),\zeta_1(x_3,x_4),\zeta_2(x_3,x_4)),
g^2(\phi_1(x_1,x_2),\phi_2(x_1,x_2),\zeta_1(x_3,x_4),\zeta_2(x_3,x_4)))
\end{eqnarray*}
(obviously, if $f^{1,2}\sim g^{1,2}$ and $f^{1,2}$ solves (\ref{Hiro}), then $g^{1,2}$ does the same).
\end{theo}

\noindent Item 1 is already argued. To prove Item 2 we let $\{\F_\la\}$ be a Kronecker web of $3$-web type. Then $T\F_{\la_1}\oplus T\F_{\la_3}=TM$ an we can find coordinates $x_1,\ldots,x_4$ such that $T\F_{\la_1}=\langle \frac{\d }{\d x_1}, \frac{\d }{\d x_2}\rangle, \F_{\la_3}=\langle \frac{\d }{\d x_3}, \frac{\d }{\d x_4}\rangle$. Define a Nijenhuis operator $N$ by $N|_{T\F_{\la_i}}=\la_i\Id_{T\F_{\la_i}}$.

It turn out that the  distributions $T\F_\la$ and $\H_\la:=(N-\la\Id_{TM})(T\F_{\infty})$ coincide for any $\la$. Indeed, both of them are of the form $\langle (\la-\la_3)X_1+(\la-\la_1)X_2,(\la-\la_3)Y_1+(\la-\la_1)Y_2\rangle$, where the vector fields $X_{1,2},Y_{1,2}$ are linearly independent everywhere. Since $\H_\la$ and $\F_\la$ coincide for $\la=\la_{1,3}$, we have $T\H_\la=\langle (\la-\la_3)X_1+(\la-\la_1)X_2,(\la-\la_3)Y_1+(\la-\la_1)Y_2\rangle$ and $T\F_\la=\langle (\la-\la_3)(a^1_1X_1+b^1_1Y_1)+(\la-\la_1)(a^1_2X_2+b^1_2Y_2),(\la-\la_3)(a^2_1X_1+b^2_1Y_1)+
(\la-\la_1)(a^2_2X_2+b^2_2Y_2)\rangle$ for some vector fields $X_{1,2},Y_{1,2}$ and functions $a^i_j,b^i_j$. On the other hand, the equality $\H_\infty=\F_\infty$, $\langle X_1+X_2,Y_1+Y_2\rangle=\langle (a^1_1X_1+b^1_1Y_1)+(a^1_2X_2+b^1_2Y_2),(a^2_1X_1+b^2_1Y_1)+
(a^2_2X_2+b^2_2Y_2)\rangle$,  implies due to the linear independence of $X_{1,2},Y_{1,2}$ that $a^1_1=a^1_2,a^2_1=a^2_2,b^1_1=b^1_2,b^2_1=b^2_2$ and $T\F_\la=\langle a^1_1[(\la-\la_3)X_1+(\la-\la_1)X_2]+b^1_1[(\la-\la_3)Y_1+(\la-\la_1)Y_2],
a^2_1[(\la-\la_3)X_1+(\la-\la_1)X_2]+b^2_1[(\la-\la_3)Y_1+(\la-\la_1)Y_2]\rangle$, which proves the claim.

In particular, $\H_0=\F_0=N(\F_\infty)$ is integrable and considerations above show that independent  first integrals of $\F_\infty$ should satisfy (\ref{Hiro}) and conditions (\ref{jac}).

The last part of the theorem can be argued in the same way as that of Theorem \ref{6}.
\qed

\smallskip

Of course the construction can be repeated for other normal forms of Nijenhuis operators in $\R^4$. However a natural question whether each Kronecker web leads to some solution of the corresponding system is more subtle (cf. Remark \ref{remn}).

%
\section{An overview of related results}
\label{notes}

Below list some related results that are beyond the scope of this paper.

For the general theory of Veronese and Kronecker webs, including their local classification see \cite{t7}, \cite{t8}, \cite{Kr1}. In the last article and in \cite{Kr2} the relations of Kronecker webs with systems of ODEs are discussed and adapted connections are built, which allow to distinguish among flat and nonflat webs (cf. the definition of flatness after Theorem \ref{turiel}).

The problem of bisymplectic realizations of generic Kronecker bihamiltonian structures   is studied in \cite{petalidou}.

B\"{a}cklund transformations, contact symmetry algebras and some exact solutions of the equations of types A--C (see Section \ref{listPDE}) and also of type D, which corresponds to the case of a Nijenhuis operator with imaginary eigenvalues and which we omit in this article, can be found in \cite{kruglikovP}.

As mentioned in Introduction in paper \cite{dunajskiKr} one can find a description of relations of Veronese webs in  3D with the hyper-CR Einstein--Weyl structures, in particular an explicit formula of such a structure based on a solution of the Hirota equation. Similar formulae for other equations of types A--D are discussed in \cite{kruglikovP}. In recent paper \cite{Kr3} a twistor geometric approach is used to treat on the same base equations of types A--D and  mentioned  in Introduction hyper-CR equation, which cannot be included in the scheme of \cite{kruglikovP}. In the same paper \cite{Kr3} there appears system (\ref{Hiro}) and its twistor geometric deformations.

Finally, these deformations, their  generalizations to higher dimensions, and relations with the Pleba\'{n}ski equation are discussed in \cite{Kr4}.

%
\section*{Appendix: Classification of cyclic Nijenhuis operators in 3D (after F.\,J.\ Turiel)}\label{app}

In papers \cite{t4,grifoneMehdi} there was obtained a local classification of complex analytic Nijenhuis (1,1)-tensors $N:TM\to TM$ (in a vicinity of a regular point \cite[p. 451]{t4}) under additional assumption of existence of a complete family of the so-called conservation laws. This assumption is equivalent to vanishing of the invariant $P_N$, which is automatically trivial in the case of cyclic $N$    \cite[p. 450]{t4}, i.e. when the space $T_xM$ is cyclic for $N_x,x\in M$. Here we recall the normal forms obtained in this case for 3-dimensional $M$.

The results of \cite{t4} imply that for any cyclic Nijenhuis (1,1)-tensor in a vicinity of a regular point $x^0$ there exist a local system of coordinates $(x_1,x_2,x_3)$ and pairwise distinct constants $a_1,a_2,a_3$ such that the coordinates $(x^0_1,x^0_2,x^0_3)$ of $x^0$ are also pairwise distinct  and  the matrix $N$ of the corresponding operator in the basis $\frac{\partial }{\partial x_1},\frac{\partial }{\partial x_2},\frac{\partial }{\partial x_3}$ is one from the following list. Besides the matrices $N$ themselves   below we list also  their Frobenius forms $F$ and  their Jordan forms $J$ .

\begin{itemize}
\item[\bf A0.] $N_{A0}=N_{A0}(x_1,x_2,x_3):=\left[
           \begin{array}{ccc}
             x_1  & 0 & 0 \\
             0 & x_2 & 0 \\
             0 & 0 & x_3 \\
           \end{array}
         \right]$,
$F_{A0}=F_{A0}(x_1,x_2,x_3)=\left[
           \begin{array}{ccc}
             0 & 0 & x_1x_2x_3 \\
             1 & 0 & -x_1x_2-x_1x_3-x_2x_3 \\
             0 & 1 & x_1+x_2+x_3 \\
           \end{array}
         \right]$,
$J_{A0}=N_{A0}$.
\item[\bf A1.] $N_{A1}:=N_{A0}(x_1,x_2,a_3),F_{A1}:=F_{A0}(x_1,x_2,a_3),J_{A1}=N_{A1}$.
\item[\bf A2.] $N_{A2}:=N_{A0}(x_1,a_2,a_3),F_{A2}:=F_{A0}(x_1,a_2,a_3),J_{A2}=N_{A2}$.
\item[\bf A3.] $N_{A3}:=N_{A0}(a_1,a_2,a_3),F_{A3}:=F_{A0}(a_1,a_2,a_3),J_{A3}=N_{A3}$.
%2
\item[\bf B0.] $N_{B0}=N_{B0}(x_2,x_3):=\left[
           \begin{array}{ccc}
              x_2  & 1 & 0 \\
             0 & x_2 & 0 \\
             0 & 0 & x_3 \\
           \end{array}
         \right]$,
$F_{B0}=F_{B0}(x_2,x_3):=\left[
           \begin{array}{ccc}
             0 & 0 & x_2^2x_3 \\
             1 & 0 & -x_2^2-2x_2x_3 \\
             0 & 1 & 2x_2+x_3 \\
           \end{array}
         \right]$,
$J_{B0}=N_{B0}$.
\item[\bf B1.] $N_{B1}:=N_{B0}(x_2,a_3),F_{B1}:=F_{B0}(x_2,a_3),J_{B1}:=N_{B1}$.
\item[\bf B2.] $N_{B2}:=N_{B0}(a_2,x_3),F_{B2}:=F_{B0}(a_2,x_3),J_{B2}:=N_{B2}$.
\item[\bf B3.] $N_{B3}:=N_{B0}(a_2,a_3),F_{B3}:=F_{B0}(a_2,a_3),J_{B3}:=N_{B3}$.
%3
\item[\bf C0.] $N_{C0}=N_{C0}(x_2,x_3):=\left[
           \begin{array}{ccc}
             x_3  & 0 & 1 \\
             1 & x_3 & -x_2 \\
             0 & 0 & x_3 \\
           \end{array}
         \right]$,
$F_{C0}=F_{C0}(x_3):=\left[
           \begin{array}{ccc}
             0 & 0 & x_3^3 \\
             1 & 0 & -3x_3^2 \\
             0 & 1 & 3x_3       \\
           \end{array}
         \right]$,
$J_{C0}=J_{C0}(x_3):=\left[
           \begin{array}{ccc}
             x_3 & 1 & 0 \\
             0 & x_3 & 1 \\
             0 & 0 & x_3 \\
           \end{array}
         \right]$.
\item[\bf C1.] $N_{C1}:=J_{C0}(a_3),F_{C1}:=F_{C0}(a_3),J_{C1}:=N_{C1}$.

\end{itemize}

%\bibliography{C:/Users/Lenovo/Documens/BIB/Main}
%\bibliography{D:/Dokumenty/BIB/Main}

\end{document}